\documentclass[10pt, a4paper]{article}

\usepackage{geometry}
\usepackage{mathtools}
\usepackage{amsmath}
\usepackage{amssymb}
\usepackage{amsfonts}
\usepackage{amsthm}
\usepackage{relsize}

\usepackage[english]{babel}
\usepackage[utf8]{inputenc}
\usepackage[T1]{fontenc}
\usepackage[hidelinks]{hyperref}
\usepackage[round]{natbib}
\usepackage[pdftex,dvipsnames,table]{xcolor}
\usepackage[nameinlink,noabbrev]{cleveref}

\usepackage{float}
\usepackage{environ}
\usepackage{tikz} 
\usepackage{tikzsymbols}

\usepackage{indentfirst}
\usepackage{url}
\usepackage{latexsym}
\usepackage{calc}
\usepackage{epsfig}
\usepackage{fancyhdr}
\usepackage{graphicx}
\usepackage{scalerel}
\usepackage{xargs}
\usepackage{xspace}
\usepackage{pdftexcmds}
\usepackage{pgfplots}
\usepackage{hyphenat}

\usepackage{etoolbox}
\usepackage{setspace}
\usepackage{bookmark}
\usepackage{booktabs}
\usepackage{multirow}

\usepackage[usestackEOL]{stackengine}

\edef\restoreparindent{\parindent=\the\parindent\relax}
\usepackage{parskip}
\restoreparindent
\usepackage{caption}
\usepackage{subcaption}

\usepackage{enumitem}

\usepackage{etoolbox}
\usepackage{lineno}
\usepackage[normalem]{ulem}
\usepackage[colorinlistoftodos,prependcaption,textsize=tiny]{todonotes}

\newcommandx{\unsure}[2][1=]{%
  \todo[linecolor=red, backgroundcolor=red!25, bordercolor=red, #1]{%
  \begin{flushleft}#2\end{flushleft}%
}}

\newcommandx{\change}[2][1=]{%
  \todo[linecolor=blue, backgroundcolor=blue!25, bordercolor=blue, #1]{%
  \begin{flushleft}#2\end{flushleft}%
}}

\newcommandx{\info}[2][1=]{%
  \todo[linecolor=OliveGreen, backgroundcolor=OliveGreen!25, bordercolor=OliveGreen, #1]{%
  \begin{flushleft}#2\end{flushleft}%
}}

\newcommandx{\improvement}[2][1=]{%
  \todo[linecolor=Plum, backgroundcolor=Plum!25, bordercolor=Plum, #1]{%
  \begin{flushleft}#2\end{flushleft}%
}}

\newcommandx{\thiswillnotshow}[2][1=]{\todo[disable, #1]{#2}}

\definecolor{RED}{rgb}{1,0,0}\definecolor{BLUE}{rgb}{0,0,1}

\makeatletter
\newif\iftx@libertine
\newif\iftx@minion
\newif\iftx@coch
\newif\iftx@ch
\newif\iftx@stxtwo
\newif\iftx@ebgm
\newif\iftx@ut
\newif\iftx@nc
\newif\iftx@ams

\newcommand{\appendixnumberline}[1]{Appendix\space}

\def\mylabel#1{\ifmeasuring@\else\ltx@label{#1}\fi}
\DeclareRobustCommand{\makeref}[1]{%
  \makephantom%
  \quad\label{#1}
}

\DeclareRobustCommand{\veb}[1]{\mathpalette\do@veb{#1}}
\newcommand{\do@veb}[2]{%
  \fix@cev{#1}{+}%
  \mbox{$\m@th#1\vec{\mbox{$\fix@cev{#1}{-}\m@th#1#2\fix@cev{#1}{+}$}}$}%
  \fix@cev{#1}{-}%
}

\DeclareRobustCommand{\cev}[1]{\mathpalette\do@cev{#1}}
\newcommand{\do@cev}[2]{%
  \fix@cev{#1}{+}%
  \reflectbox{$\m@th#1\vec{\reflectbox{$\fix@cev{#1}{-}\m@th#1#2\fix@cev{#1}{+}$}}$}%
  \fix@cev{#1}{-}%
}

\newcommand{\fix@cev}[2]{%
  \ifx#1\displaystyle 
    \mkern#23mu
  \else
    \ifx#1\textstyle
      \mkern#23mu
    \else
      \ifx#1\scriptstyle
        \mkern#22mu
      \else
        \mkern#22mu
      \fi
    \fi
  \fi
}

\makeatother

\newcommand{\RemoveSpaces}[1]{%
  \begingroup%
  \spaceskip=1.2mm%
  \xspaceskip=1.2mm%
  \mbox{#1}%
  \endgroup
}

\DeclareMathAlphabet{\matheuler}{U}{zeur}{m}{n}
\DeclareSymbolFont{lettersslanted}{OML}{zplm}{m}{n}
\DeclareMathSymbol{p}{\mathalpha}{lettersslanted}{`p}

\newcommand{\shortname}{EVSP\xspace}
\newcommand{\ie}{i.e.\xspace}

\newcommand{\incoming}[1]{\!\smash{\cev{\,#1}}}
\newcommand{\outgoing}[1]{\!\!\smash{\veb{\,\,#1}}}

\newcommand{\basephantom}{\incoming{\matheuler{E}}}
\newcommand{\makephantom}{\vphantom{\sum\limits_{\mathrlap{\substack{\basephantom}}}}}
\newcommand{\centerlap}[1]{\limits_{\mathclap{\substack{#1\vphantom{\basephantom}}}}}
\newcommand{\rightlap}[1]{\limits_{\mathrlap{\substack{\mathsurround=-7pt#1\vphantom{\basephantom}}}}}

\newcommand{\lowprime}{\mkern0mu\raise0.6ex\hbox{$\scriptstyle\prime$}}
\newcommand{\CC}{%
  \sffamily C\nolinebreak\hspace{-.01em}\raisebox{.3ex}{\footnotesize +}%
             \nolinebreak\hspace{-.01em}\raisebox{.3ex}{\footnotesize +}%
}

\newcommand{\eqrefsmallsize}[1]{\eqref{#1}}

\newcounter{fc}
\newcounter{sc}

\renewcommand{\thefc}{\shortname\!\!\arabic{fc}}
\renewcommand{\thesc}{\shortname\!\!\arabic{sc}-S}

\newcommand{\makerefF}[1]{\vspace*{0.15cm}\refstepcounter{fc}%
  \noindent(\hspace*{.2mm}\thefc\label{#1})\xspace
}
\newcommand{\makerefS}[1]{\vspace*{0.15cm}\refstepcounter{sc}%
  \noindent(\hspace*{.2mm}\thesc\label{#1})\xspace
}
\newcommand{\makerefSb}[1]{\refstepcounter{sc}%
  \thesc\label{#1}\xspace
}

\newcommand{\spacetimezeroskips}[2]{\hspace*{-.3mm}(\!#1_{\hspace*{-.3mm}#2}\!)}
\newcommand{\doublespacetimezeroskips}[4]{%
  (\hspace*{-.3mm}#1_{\hspace*{-.3mm}#2}\!,\hspace*{-.1mm}#3_{\hspace*{-.3mm}#4}\!)%
}

\newcommand{\cardinality}[1]{\vert\parbox{.8em}{\centering$#1$}\vert}

\newcolumntype{C}[1]{>{\centering\arraybackslash}m{#1}}

\newcommand{\cStack}[1]{%
  \Centerstack[l]{#1}%
  \vspace*{.02cm}%
}
\newcommand{\mStack}[1]{%
\begin{minipage}{.3\textwidth}%
  \raggedright%
  #1%
\end{minipage}\vspace*{2mm}%
}

\newcommand{\rotateBox}[1]{%
  \vspace*{-.2cm}%
  \rotatebox{40}{\Centerstack[l]{\hspace*{-.35cm}#1}}\hspace*{-.45cm}%
  \vspace*{-.1cm}%
}

\newtheorem{corollary}{Corollary}
\newtheorem{proposition}{Proposition}

\pgfplotsset{
  box plot/.style={
    /pgfplots/.cd,
    black,
    only marks,
    mark=-,
    mark size=\pgfkeysvalueof{/pgfplots/box plot width},
    /pgfplots/error bars/y dir=plus,
    /pgfplots/error bars/y explicit,
    /pgfplots/table/x index=\pgfkeysvalueof{/pgfplots/box plot x index},
  },
  box plot box/.style={
    /pgfplots/error bars/draw error bar/.code 2 args={%
      \draw  ##1 -- ++(\pgfkeysvalueof{/pgfplots/box plot width},0pt) %
         |- ##2 -- ++(-\pgfkeysvalueof{/pgfplots/box plot width},0pt) |- ##1 -- cycle;
    },
    /pgfplots/table/.cd,
    y index=\pgfkeysvalueof{/pgfplots/box plot box top index},
    y error expr={
      \thisrowno{\pgfkeysvalueof{/pgfplots/box plot box bottom index}}
      - \thisrowno{\pgfkeysvalueof{/pgfplots/box plot box top index}}
    },
    /pgfplots/box plot
  },
  box plot top whisker/.style={
    /pgfplots/error bars/draw error bar/.code 2 args={%
      \pgfkeysgetvalue{/pgfplots/error bars/error mark}%
      {\pgfplotserrorbarsmark}%
      \pgfkeysgetvalue{/pgfplots/error bars/error mark options}%
      {\pgfplotserrorbarsmarkopts}%
      \path ##1 -- ##2;
    },
    /pgfplots/table/.cd,
    y index=\pgfkeysvalueof{/pgfplots/box plot whisker top index},
    y error expr={
      \thisrowno{\pgfkeysvalueof{/pgfplots/box plot box top index}}
      - \thisrowno{\pgfkeysvalueof{/pgfplots/box plot whisker top index}}
    },
    /pgfplots/box plot
  },
  box plot bottom whisker/.style={
    /pgfplots/error bars/draw error bar/.code 2 args={%
      \pgfkeysgetvalue{/pgfplots/error bars/error mark}%
      {\pgfplotserrorbarsmark}%
      \pgfkeysgetvalue{/pgfplots/error bars/error mark options}%
      {\pgfplotserrorbarsmarkopts}%
      \path ##1 -- ##2;
    },
    /pgfplots/table/.cd,
    y index=\pgfkeysvalueof{/pgfplots/box plot whisker bottom index},
    y error expr={
      \thisrowno{\pgfkeysvalueof{/pgfplots/box plot box bottom index}}
      - \thisrowno{\pgfkeysvalueof{/pgfplots/box plot whisker bottom index}}
    },
    /pgfplots/box plot
  },
  box plot median/.style={
    /pgfplots/box plot,
    /pgfplots/table/y index=\pgfkeysvalueof{/pgfplots/box plot median index}
  },
  box plot width/.initial=1em,
  box plot x index/.initial=0,
  box plot median index/.initial=1,
  box plot box top index/.initial=2,
  box plot box bottom index/.initial=3,
  box plot whisker top index/.initial=4,
  box plot whisker bottom index/.initial=5,
  x tick label style={/pgf/number format/.cd,fixed,precision=3, set thousands separator={}}
}

\usetikzlibrary{spy,positioning,arrows,arrows.meta}
  
\newcommand{\boxplot}[2][]{
  \addplot [box plot median,#1] table {#2};
  \addplot [forget plot, box plot box,#1] table {#2};
  \addplot [forget plot, box plot top whisker,#1] table {#2};
  \addplot [forget plot, box plot bottom whisker,#1] table {#2};
}

\title{\LARGE%
 On the complexity and modeling of the electric vehicle sharing problem
}

\author{Welverton R. Silva, Fábio L. Usberti, Rafael C.S. Schouery\thanks{Corresponding author.%
\newline\hspace*{.5cm}\textit{E-mail address}: \href{rafael@ic.unicamp.br}{rafael@ic.unicamp.br}
}}

\begin{document}

\date{}
\maketitle

\vspace{-20pt}
\begin{center}\footnotesize 
  Institute of Computing, University of Campinas\\Campinas, São Paulo, Brazil
\end{center}

\bigskip\bigskip\noindent
{\small {\bf ABSTRACT.}
We introduce the electric vehicle sharing problem~(\shortname), a problem that arises from the planning and operation of electric car-sharing systems which allow one-way rental of vehicles.~The problem aims at finding the maximum total daily rental time in which customers' demands are assigned to the existing fleet.~In addition, either all of the customer's demands are completely fulfilled or the customer does not use the system at all.~We show that the~\shortname is~\mbox{NP-hard}, and we provide four mixed-integer linear programming formulations based on space-time network flow models, along with some theoretical results.~We perform a comprehensive computational study of the behavior of the proposed formulations using two benchmark sets, one of which is based on real-world data from an electric car-sharing system located in~Fortaleza, Brazil.~The results show that our best formulation is effective in solving instances where each customer has only one demand.~In general, we were able to optimally solve at least~\(55\%\) of the instances within the time limit of one hour.}


\baselineskip=\normalbaselineskip
\sloppy

\section{Introduction}
\label{sec:introduction}

In car-sharing systems, cars are used on demand by customers for short periods to travel relatively small distances~(short-term car rental).~Usually, cars are spread throughout the city in parking stations or on-street parking places, and customers can book a car through a smartphone app~(subject to availability).~Car-sharing services provide a low-cost and highly flexible alternative to private mobility, freeing people from the costs and responsibilities associated with owning a car.~These services also improve the quality of urban transport systems by reducing congestion, greenhouse gas emissions, parking shortages, and noise pollution~\citep{Shaheen/2007,Zakaria/2014}.~An example is the City of~Bremen in northwestern~Germany, which has been successfully integrating car-sharing systems into urban and transport planning for more than a decade to reduce the number of cars circulating and the pressure on the limited number of parking spaces~\citep{Glotz/2016}.

Car-sharing services can be~\emph{free-floating} or~\emph{station-based} systems.~In free-floating systems, cars can be parked in free parking spaces, and customers can pick up available cars and drop them off in any free parking space within a given zone or covered area.~In station-based systems, cars are allocated to dedicated stations for parking.~A station-based system can be classified into two types of service models, namely~\emph{one-way} and~\emph{two-way}.~One-way rental is more flexible for customers as the rented car can be returned to a station of your choice.~Two-way rental, on the other hand, is more restricted as customers have to return the cars to their original pick-up stations~\citep{Nourinejad/2015:May}.

In this paper, we focus on one-way station-based electric car-sharing systems, motivated by their establishment in many cities around the world~\citep{Hall/2017}.~For example, the Alternative Vehicles for~Mobility~(VAMO, acronym for\,\emph{Veículos Alternativos para Mobilidade}) is the first public one-way electric car-sharing system in Brazil, located in the city of~Fortaleza since~2016.~VAMO Fortaleza is an initiative aiming to promote sustainable urban mobility, partly in response to the excessive use of private cars in Fortaleza, which increased urban problems such as noise and pollution~\citep{Pereira/2021}.~After one year of operation,~\cite{Teles/2018} estimated that this initiative reduced CO\(_2\) emissions to the atmosphere by about~\(5.5\) tons.

Although instant access (\ie, with no advanced reservation) offers flexibility for customers, it decreases the ability of car-sharing service providers to predict demands~\citep{Brendel/2017}.~Thus, we consider reservation-based car-sharing systems, that is, systems in which customers must book cars before using them, stating the pick-up and drop-off stations of their trips.

This work also considers customers with multiple driving demands.~For example, a customer with two demands could have one demand from a station in a residential area to a station in a central business district and, after a working day, return to the residential area with another driving demand.~In this case, a system that does not guarantee that both demands will be fulfilled could be detrimental to the customer's needs.~Also, for the car-sharing service provider, the planning of trips (including routes and timetables of the cars) can be severely compromised if customers should cancel their trips as a result of some of their demands not being fulfilled.~Thus, we consider a system where either all or neither of the demands of a customer are satisfied.

We introduce a new combinatorial optimization problem, the~electric vehicle sharing \mbox{problem}~(\shortname), which implements a staffless\footnote{Staffless refers to a car-sharing system that does not require dedicated staff (operators) moving the cars for charging, repositioning, and rebalancing purposes.} station-based car-sharing system that allows one-way rental of vehicles.~In the~\shortname, the driving demands for electric vehicles~(EVs) are known in advance (say, the day before), and the objective is to maximize the total daily vehicle rental time by fulfilling customers' demands.~We assume that the customer can have more than one demand and is not interested in partial rentals, that is, either all customer's demands are completely fulfilled or the customer does not use the system at all.~If~all the customer's demands are fulfilled, we say that the customer is~\emph{served}.

In addition to the introduction of a new problem arising in car-sharing systems, our main contributions in this paper can be summarized as follows.

\begin{enumerate}[label=({\roman*})]
  \item Four mixed-integer linear programming formulations for the~\shortname based on homogeneous and heterogeneous space-time networks, and theoretical results for the dominance relation of the linear relaxation of the formulations.
  \item An \mbox{NP-hardness} proof for the~\shortname.
  \item A benchmark of randomly generated instances and another benchmark based on data from the~VAMO Fortaleza system~\citep{Fortaleza/2018}.
  \item A comprehensive computational study on our instances to compare the performance of our formulations using a commercial solver for mixed-integer linear programming.
\end{enumerate}

The~\shortname distinguishes itself from previous vehicle sharing problems by considering the energy consumption and recharging of EV batteries, which significantly increases the complexity of developing efficient and robust exact algorithms.~In addition, the~\shortname is a problem whose interest is not limited to car-sharing systems.~Other applications of the~\shortname include the rental of micromobility\footnote{Micromobility refers to small and lightweight vehicles to travel short distances around cities, most of which are used individually, such as the use of bicycles and scooters~\citep{Sengul/2021}.} electric vehicles, such as docked e-bikes and electric scooters.~Also, the~\shortname can be used to model drone one-way path planning problems using recharging stations for pick-up and delivery services~\citep{Huang/2021, Pachayappan/2021}.


The remainder of this paper is organized as follows.~In~\Cref{sec:literature}, we survey the literature and related works.~In~\Cref{sec:problem_definition}, we give a formal definition of the~\shortname.~In~\Cref{sec:formulations}, we propose mixed-integer linear programming formulations to model the problem.~In \Cref{sec:complexity}, we~present the \mbox{NP-hardness} proof for the~\shortname.~In~\Cref{sec:experiments}, we perform a detailed computational study.~Finally, in~\Cref{sec:conclusions}, we give concluding remarks and an outlook on future research.

\section{Previous and related work}
\label{sec:literature}

We dedicate this section to briefly describe the problem of which ours is a variant, and survey related work in the literature focusing on one-way electric car-sharing systems.

\citet{Bohmova/2016} propose the~transfers for commuting~(TFC) problem, which is motivated by applications such as car-sharing with the objective of maximizing the number of customers served in a one-way rental system without relocations.~There are two stations~\(\textsf{A}\) and~\(\textsf{B}\), and each customer has \mbox{exactly} two demands~(rental requests) in opposite directions, that is, one from station~\(\textsf{A}\) to~\(\textsf{B}\) and the other from station~\(\textsf{B}\) to~\(\textsf{A}\), although not necessarily in this order.~The following assumptions are made: the start and end times of all demands are known in advance, a fleet of homogeneous vehicles is initially distributed between both stations, and the number of parking spaces at each station is sufficient to cover all vehicles in the system.

\begin{figure}[H]
  \centering
  \resizebox{\textwidth}{!}{%
    \tikzset{every picture/.style={line width=0.75pt}} 

\begin{tikzpicture}[x=0.75pt,y=0.75pt,yscale=-1,xscale=1]

   \draw [draw opacity=0][fill={rgb, 255:red, 255; green, 210; blue, 210}, fill opacity=1] (24.17,145.33) -- (557.85,145.33) -- (557.85,135.33) -- (580.56,155.33) -- (557.85,175.33) -- (557.85,165.33) -- (24.17,165.33) -- cycle;
   \draw [draw opacity=0][fill={rgb, 255:red, 185; green, 215; blue, 252}, fill opacity=1] (24.67,43.33) -- (557.65,43.33) -- (557.65,33.33) -- (580.33,53.33) -- (557.65,73.33) -- (557.65,63.33) -- (24.67,63.33) -- cycle;

   \draw [->,>={Triangle[scale=.6]},line width=.35mm,color={rgb, 255:red, 245; green, 166; blue, 35}, draw opacity=1] (56.67,59.33) -- (125.83,149.5);
   \draw [->,>={Triangle[scale=.6]},line width=.35mm,color={rgb, 255:red, 245; green, 166; blue, 35}, draw opacity=1] (309,149.75) -- (414.71,59);
   \draw [->,>={Triangle[scale=.6]},line width=.35mm,color={rgb, 255:red, 145; green, 45; blue, 123}, draw opacity=1] (158.33,149.75) -- (198.74,59);
   \draw [->,>={Triangle[scale=.6]},line width=.35mm,color={rgb, 255:red, 145; green, 25; blue, 123}, draw opacity=1] (305.5,59.33) -- (375.12,149.5);
   \draw [->,>={Triangle[scale=.6]},line width=.35mm,color={rgb, 255:red, 36; green, 114; blue, 18}, draw opacity=1] (203,150.17) -- (275.59,59);
   \draw [->,>={Triangle[scale=.6]},line width=.35mm,color={rgb, 255:red, 36; green, 114; blue, 18}, draw opacity=1] (477.87,59.33) -- (538.8,149.5);
   \draw [->,>={Triangle[scale=.6]},line width=.35mm,color={rgb, 255:red, 179; green, 35; blue, 24}, draw opacity=1] (225.5,59.33) -- (280.92,149.5);
   \draw [->,>={Triangle[scale=.6]},line width=.35mm,color={rgb, 255:red, 179; green, 35; blue, 24}, draw opacity=1] (421.67,149.75) -- (451.05,59);
      
   \draw (50,45) node [anchor=north west][inner sep=0.75pt] [align=left] {$\displaystyle t_{1}$};
   \draw (120,150) node [anchor=north west][inner sep=0.75pt] [align=left] {$\displaystyle t_{2}$};
   \draw (149,150) node [anchor=north west][inner sep=0.75pt] [align=left] {$\displaystyle t_{3}$};
   \draw (196,45) node [anchor=north west][inner sep=0.75pt] [align=left] {$\displaystyle t_{4}$};
   \draw (196,150) node [anchor=north west][inner sep=0.75pt] [align=left] {$\displaystyle t_{4}$};
   \draw (218.18,45) node [anchor=north west][inner sep=0.75pt] [align=left] {$\displaystyle t_{5}$};
   \draw (273,150) node [anchor=north west][inner sep=0.75pt] [align=left] {$\displaystyle t_{6}$};
   \draw (273,45) node [anchor=north west][inner sep=0.75pt] [align=left] {$\displaystyle t_{6}$};
   \draw (300,45) node [anchor=north west][inner sep=0.75pt] [align=left] {$\displaystyle t_{7}$};
   \draw (300,150) node [anchor=north west][inner sep=0.75pt] [align=left] {$\displaystyle t_{7}$};
   \draw (369,150) node [anchor=north west][inner sep=0.75pt] [align=left] {$\displaystyle t_{8}$};
   \draw (414,45) node [anchor=north west][inner sep=0.75pt] [align=left] {$\displaystyle t_{9}$};
   \draw (414,150) node [anchor=north west][inner sep=0.75pt] [align=left] {$\displaystyle t_{9}$};
   \draw (444.51,45) node [anchor=north west][inner sep=0.75pt] [align=left] {$\displaystyle t_{10}$};
   \draw (470.51,45) node [anchor=north west][inner sep=0.75pt] [align=left] {$\displaystyle t_{11}$};
   \draw (528.68,150) node [anchor=north west][inner sep=0.75pt] [align=left] {$\displaystyle t_{12}$};

   \draw (557.85,45) node [anchor=north west][inner sep=0.75pt] [align=left] {$\textsf{A}$};
   \draw (557.85,150) node [anchor=north west][inner sep=0.75pt] [align=left] {$\textsf{B}$};

\end{tikzpicture}
  }
  \caption{An example of the~TFC problem.~Each station is represented by a horizontal arrow showing the progression of time, discretized as a set~\mbox{\(\{t_{1}, \dots, t_{12}\}\)} of time instants corresponding to the start and end times of the customers' demands.~The demands are displayed as full arrows between the two stations, and all demands with the same color represent the same customer.}
  \label{fig:TFC-example}
\end{figure}
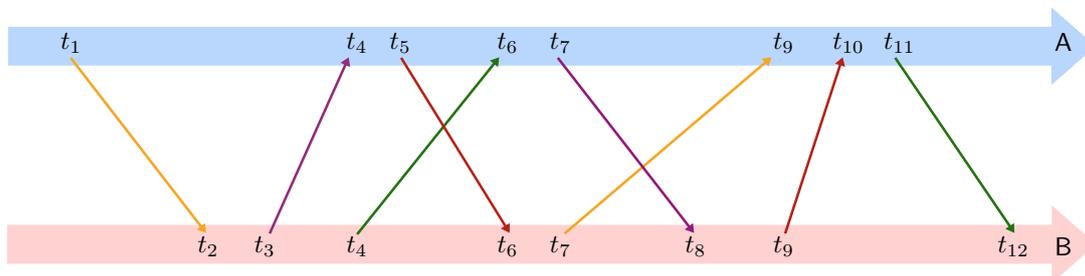

The~TFC problem aims at maximizing the number of customers with demands fulfilled by a fleet of vehicles, such that each customer either has both demands fulfilled or does not use the system at~all.~In addition, the demands of a customer can be fulfilled by different vehicles.~To illustrate this problem,~\Cref{fig:TFC-example} gives an example with four customers.~In this example, each demand is modeled as a vehicle commute between the two stations, that is, a driving demand represents a request for a vehicle to move from one station to another station at the given time interval.

\citet{Bohmova/2016} prove that the~TFC problem is~\mbox{NP-hard}, and also~\mbox{APX-hard}, even if all demands take the same time to commute between two stations, and there is only one vehicle.~The proof is based on a reduction from a variant of~\mbox{MAX-3SAT}~\citep{Ausiello/2012}.~In contrast, the variant of the~TFC problem considered in their work, where each customer has only one demand, is solvable in polynomial time by reducing to the minimum-cost maximum-flow problem.



\citet{Silva/2020} present a mixed-integer linear programming formulation for the~TFC problem.~To this end, the authors model the problem as a network flow problem and present a preprocessing procedure to reduce the number of nodes and arcs in the network, thereby significantly reducing the number of decision variables and constraints in the integer linear program.~Computational experiments are conducted with randomly generated instances to evaluate the quality of the proposed approaches.~The results show that the proposed formulation produces solutions with small optimality gaps for large instances.


Recently, the~TFC problem was formally modeled as an online problem by~\citet{Luo/2020}.~In this online version, each customer has a pair of demands in opposite directions, where the decision on each pair must be made without knowledge of future demands immediately at the time when the pair is submitted.~For this online version of the~TFC problem, the authors propose two algorithms, one of which is \mbox{a \(4\)-competitive} algorithm.

As can be seen in \Cref{tab:studies}, most of the existing one-way car-sharing works focus on the problem caused by the imbalanced spatial and temporal distribution of vehicles among stations.~In this case, the regular relocation of vehicles between stations is necessary to ensure that there are adequate amounts of vehicles available at stations to serve customers' demands~\citep{Kek/2009}.~Vehicle relocation can be carried out by operators or by the customers themselves by following instructions on where to park the vehicle through a price incentive mechanism aimed at maximizing profit or minimizing cost.~These strategies are commonly known as~\emph{operator-based} and~\emph{user-based} vehicle relocation, respectively.~A literature review on the vehicle relocation problem in one-way car sharing is presented in~\mbox{\citet{Illgen/2019}}.

\strutlongstacks{T}
\begin{table}[H]
  \singlespacing
  \caption{Summary of recent studies on one-way station-based electric car-sharing systems.}
  \label{tab:studies}
  \vspace*{.5mm}
  \resizebox{\textwidth}{!}{%
    \begin{tabular}{lllC{1.3cm}C{1.5cm}C{1cm}C{2.8cm}c}
    \toprule\rule{0pt}{2\normalbaselineskip}
    \cStack{Research article}       & \cStack{Focus}                                               & \cStack{Objectives}                                                        & \rotateBox{Multiple driving\\demands} & \rotateBox{~~~~Reservation-based~~~~~\\~~~~~~system} & \rotateBox{Staffless station}\hspace*{-10mm} & \rotateBox{~~Partially charged\\\hspace*{.39cm}vehicle rentals} \vspace*{-4mm} \\ \hline\rule{0pt}{2\normalbaselineskip}
    \!Current study                 & \mStack{EV fleet operation planning; non-partial rentals}    & \mStack{Maximize total rental time}                                        & \(\bullet\)                           & \(\bullet\)                                          & \(\bullet\)                                  & \(\bullet\)                                                                    \\ 
    \cite{Ait-Ouahmed/2018}         & \mStack{Operator-based relocation; quantity of resources}    & \mStack{Maximize the number of served customers}                           &                                       & \(\bullet\)                                          &                                              &                                                                                \\
    \cite{Boyaci/2017}              & \mStack{Operator-based relocation; staff operation planning} & \mStack{Maximize number of fulfilled demands and minimize relocation cost} &                                       & \(\bullet\)                                          &                                              &                                                                                \\
    \cite{Brandstatter/2017}        & \mStack{Location analysis; charging station location}        & \mStack{Maximize the expected profit}                                      &                                       &                                                      & \(\bullet\)                                  &                                                                                \\
    \cite{Bruglieri/2017}           & \mStack{Operator-based relocation; routes and schedules}     & \mStack{Maximize the total profit}                                         &                                       & \(\bullet\)                                          &                                              & \(\bullet\)                                                                    \\
    \cite{Calik/2019}               & \mStack{Location analysis; charging station location}        & \mStack{Maximize the expected profit}                                      &                                       &                                                      & \(\bullet\)                                  &                                                                                \\ 
    \cite{Gambella/2018}            & \mStack{Operator-based relocation}                           & \mStack{Maximize the profit associated with fulfilled demands}             &                                       & \(\bullet\)                                          &                                              & \(\bullet\)                                                                    \\
    \cite{Huang/2020}               & \mStack{Operator-based relocation; user-based relocation}    & \mStack{Maximize the total profit}                                         &                                       &                                                      &                                              & \(\bullet\)                                                                    \\
    \cite{Pantelidis/2021}          & \mStack{Operator-based relocation; rebalancing policy}       & \mStack{Minimize the cumulative operating cost}                            &                                       &                                                      &                                              & \(\bullet\)                                                                    \\
    \cite{Xu/2019}                  & \mStack{Fleet size; vehicle relocation operations}           & \mStack{Maximize the profit of operators}                                  &                                       &                                                      &                                              & \(\bullet\)                                                                    \\
    \cite{Zhang/2019}               & \mStack{EV assignment with relays}                           & \mStack{Maximize the total proﬁt over the given time horizon}              &                                       & \(\bullet\)                                          &                                              & \(\bullet\)                                                                    \\
    \cite{Zhao/2018}                & \mStack{Operator-based relocation; allocation plan}          & \mStack{Minimize the total operation cost}                                 &                                       & \(\bullet\)                                          &                                              & \(\bullet\)                                                                    \\ \toprule
    \end{tabular}}
\end{table}

To the best of our knowledge, the~\shortname is the first one-way station-based electric car-sharing problem that addresses multiple driving demands.~This is a relevant generalization because it increases customer satisfaction and provides predictability for car-sharing service providers.~Also, the~\shortname  distinguishes itself from previous works on electric car-sharing by considering staffless station-based systems, which have the benefit of a reduced operational cost by not requiring frequent vehicle relocation~\citep{Gambella/2018, Folkestad/2020}.

\section{Problem definition}
\label{sec:problem_definition}

In this section, we present a formal definition of the~\shortname.~To this end, in addition to the assumptions about customers~(presented in~\Cref{sec:introduction}), we make the following assumptions about the shared electric vehicle system.

\begin{enumerate}[label=({\roman*})]
  \item The customer picks an EV at the demand start time and delivers it immediately before the demand end time.~Thus, the EV becomes available to other customers at the demand end time.
  \item Every station has a limited number of parking spaces, which some or all are equipped with charging facilities.
  \item There are no standby operators at each station, which means that the returning customer plugs in the charging cable to the~EV (if there is one at the parking space) and, the~EV is kept plugged in until it is used by the next customer.
  \item All the EVs are identical and are in the same working condition.~Also, all parking spaces are equipped with the same charging facilities, where~\mbox{\(\mu \in \mathbb{Q}_{>0}\)} denotes the energy supplied by a charging facility per unit time~(in watts per minute).
\end{enumerate}


Henceforth, consider the superscript~\(\outgoing{~}\)\, as an ``outgoing from'' labeling and the superscript~\(\incoming{~}\)\, as an ``incoming into'' labeling to indicate that a vehicle is leaving or arriving at a station, respectively.

Given these assumptions and some preliminary notations that we will use throughout this paper, the~\shortname is formally defined as follows.

Let~\mbox{\(\matheuler{T} = \{t_{1}, \dots, t_{m}\}\)} be a set of discrete time instants~(in minutes), not necessarily uniformly distributed, over a planning horizon~limited by the duration of a working day.~Let~\(\matheuler{C}\) and~\(\matheuler{S}\) denote the sets of customers and stations, respectively.

For each customer~\mbox{\(c \in \matheuler{C}\)}, let~\(\matheuler{D}_{c}\) be the set of demands for the customer \(c\).~Each demand in~\(\matheuler{D}_{c}\) is described by a quintuple~\mbox{\((\outgoing{s}, t_{i}, \incoming{s}, t_{j}, \varepsilon) \in (\matheuler{S} \times \matheuler{T})^2 \times \mathbb{Q}_{>0}\)}, such that~\mbox{\(t_{i} < t_{j}\)}, where~\(\outgoing{s}\)~stands for the pick-up station,~\(t_{i}\) indicates the departure time from the pick-up station, \(\incoming{s}\) stands for the drop-off station,~\(t_{j}\) indicates the arrival time at the drop-off station, and~\(\varepsilon\) is the energy required on the rental period, that is, an estimate of battery consumption.~There are no overlapping rental periods in the same set of demands.~Let~\mbox{\(\matheuler{D} = \bigcup_{c \in \matheuler{C}} \matheuler{D}_{c}\)} denote the set of all demands.

For each station~\mbox{\(s \in \matheuler{S}\)}, let~\(\matheuler{V}_{s}\) be the set of~EVs initially located at~\(s\) --- each one with battery capacity~\(\textsf{L}\)~(in watt-minutes).~Also, each station~\mbox{\(s \in \matheuler{S}\)} has a capacity~\(\textsf{C}_{s}\)~(\ie, the number of parking spaces), and a number~\(\textsf{R}_{s}\) of charging facilities, such that~\mbox{\(\textsf{R}_{s} \leq \textsf{C}_{s}\)} and~\mbox{\(|\matheuler{V}_{s}| \leq \textsf{C}_{s}\)}.~Let~\mbox{\(\matheuler{V} = \bigcup_{s \in \matheuler{S}} \matheuler{V}_{s}\)} denote the set of all~EVs.~Hereafter, the EVs are simply \mbox{referred} to as vehicles.

The objective of the~\shortname is to maximize the total daily vehicle rental time by assigning a set of driving demands to each vehicle in the existing fleet, such that, whenever a demand of a~customer is fulfilled, all demands of the customer must also be fulfilled.

A set of driving demands is associated with a feasible~\emph{assignment plan} for a vehicle~\mbox{\(v \in \matheuler{V}\)} if it satisfies the following conditions.~For each demand~\mbox{\(d = (\outgoing{s}, t_{i}, \incoming{s}, t_{j}, \varepsilon)\)} of a given set of driving demands,~\(d\) is fulfilled if: the vehicle~\(v\) is at the pick-up station~\(\outgoing{s}\) at time~\(t_{i}\); the battery energy of~\(v\) is greater than or equal to the required energy~\(\varepsilon\); and there is a parking space at the drop-off station that is available immediately before time~\(t_{j}\).~An assignment plan must determine whether and when the vehicle~\(v\) will be plugged into a charging facility if it is available.

As every demand starts at the beginning and ends at immediately before the end of the rental period, when a demand~\mbox{\((\outgoing{s}, t_{i}, \incoming{s}, t_{j}, \varepsilon) \in \matheuler{D}\)} is fulfilled by a vehicle~\mbox{\(v \in \matheuler{V}\)}, then the vehicle~\(v\) will be available at station~\(\incoming{s}\) for the next customer starting at time~\(t_{j}\).


\section{Mathematical formulations}
\label{sec:formulations}

We now provide mathematical formulations for the~\shortname with the following preliminary notations.

For each station~\mbox{\(s \in \matheuler{S}\)}, let~\mbox{\(\matheuler{V}^{\prime}_{s} \subseteq \matheuler{V}_{s}\)} denote the set of vehicles at station~\(s\) that are initially in parking spaces equipped with charging facilities.~Also, for each vehicle~\mbox{\(v \in \matheuler{V}\)}, let~\mbox{\(\textsf{L}_{v}\)} denote the initial energy stored in the battery of the vehicle~\(v\) at the beginning of the planning horizon.

Moreover, let~\mbox{\(\matheuler{I} = \{0, \dots, m\}\)} denote the set of indices for the time instants, where time~\(t_{0}\) is used to indicate the initial configuration of the system.~When the required energy~\(\varepsilon\) of the demand is not necessary, we omit it from the notation for brevity and denote the demand as a quadruple~\mbox{\((\outgoing{s}, t_{i}, \incoming{s}, t_{j} ) \in \matheuler{D}\)}.

\subsection{Homogeneous space-time networks model}

The~\shortname is first formulated using a flow formulation on a space-time network, in which each flow in the network represents an assignment plan for a particular vehicle.

\subsubsection*{Network representation and notation}

A network for this formulation is defined by a collection of node sets~\(\matheuler{N}_{s}\), for all~\mbox{\(s \in \matheuler{S}\)}, and a collection of arc sets~\(\matheuler{E}_{c}\) and~\(\matheuler{E}_{s}\), for all~\mbox{\(c \in \matheuler{C}\)} and~\mbox{\(s \in \matheuler{S}\)}.~In each~\(\matheuler{N}_{s}\), there is a~\mbox{\emph{space-time}} \emph{node}~\(s_{i}\), for all~\mbox{\(i \in \matheuler{I}\)}.~For simplicity and in the absence of ambiguity, we write~\(s_{i}\)~without specifying the set to which it belongs.~For each demand \mbox{\((\outgoing{s}, t_{i}, \incoming{s}, t_{j}, \varepsilon) \in \matheuler{D}_{c}\)}, there is a~\emph{demand arc} from node~\(\outgoing{s}_{i}\) to~\(\incoming{s}_{j}\) in each~\(\matheuler{E}_{c}\), where~\mbox{\(e = (\outgoing{s}_{i}, \incoming{s}_{j})\)} has the associated amount of energy consumed by the rental, denoted by~\mbox{\(\textsf{E}_{e} = \varepsilon\)}.~Finally, there is a~\emph{connecting arc} from node~\(s_{i - 1}\) to~\(s_{i}\) in each~\(\matheuler{E}_{s}\), where~\mbox{\(e = (s_{i - 1}, s_{i})\)} has the associated amount of energy recharged at station~\(s\), denoted by~\mbox{\(\textsf{E}_{e} = \mu(t_{i}-t_{i - 1})\)}, for all~\mbox{\(i \in \matheuler{I} \!\setminus\! \{0\}\)}.

To illustrate the representation of the space-time network, \Cref{fig:network} gives an example.~As the trajectory of the vehicle in the rental period is irrelevant for the problem, we are modeling the renting as a pick-up and delivery request at the given time interval.~Note that there are demands in which vehicles leave and arrive at the same station (\ie, vehicles move only in time).

\begin{figure}[H]
  \centering
  \resizebox{\textwidth}{!}{%
    \tikzset{every picture/.style={line width=0.75pt}}

\begin{tikzpicture}[x=0.75pt, y=0.75pt, yscale=-1, xscale=0.96]

   \draw [->,>={Triangle[scale=.6]},line width=.35mm,color={rgb, 255:red, 209; green, 53; blue, 53}, draw opacity=1] (133.5,62.63) -- (193.38,148.75);
   \draw [->,>={Triangle[scale=.6]},line width=.35mm,color=CadetBlue, draw opacity=1] (130.4,65.42) -- (188.81,149.89);
   \draw [->,>={Triangle[scale=.6]},line width=.35mm,color={rgb, 255:red, 70; green, 155; blue, 36}, draw opacity=1] (226.4,64.8) -- (313.89,149.5);
   \draw [->,>={Triangle[scale=.6]},line width=.35mm,color=CadetBlue, draw opacity=1] (223.9,150.73) -- (236.56,116.47) -- (254.19,66.73);
   \draw [->,>={Triangle[scale=.6]},line width=.35mm,color=BurntOrange, draw opacity=1] (100.07,158.) .. controls (154.18,172.92) and (202.3,173.42) .. (252.92,159.56);
   \draw [->,>={Triangle[scale=.6]},line width=.35mm,color={rgb, 255:red, 145; green, 25; blue, 123}, draw opacity=1] (410.07,65.) -- (561.6,149.29);
   \draw [->,>={Triangle[scale=.6]},line width=.35mm,color=BurntOrange, draw opacity=1] (347.73,150.9) -- (498.1,66.89);
   \draw [->,>={Triangle[scale=.6]},line width=.35mm,color={rgb, 255:red, 70; green, 155; blue, 36}, draw opacity=1] (410.07,150.9) -- (502.91,67.57);
   \draw [->,>={Triangle[scale=.6]},line width=.35mm,color={rgb, 255:red, 209; green, 53; blue, 53}, draw opacity=1] (595.07,151.07) -- (625.52,67.06);
   \draw [->,>={Triangle[scale=.6]},line width=.35mm,color=CadetBlue, draw opacity=1] (316.07,57.57) .. controls (376.58,37.55) and (422.12,39.2) .. (469.09,57.23);

   \draw [color={rgb, 255:red, 74; green, 144; blue, 226}, draw opacity=1][line width=1.5] (51.5,185.25) -- (664.5,185.25) (69,36) -- (69,199.2) (657.5,180.25) -- (664.5,185.25) -- (657.5,190.25) (64,43) -- (69,36) -- (74,43) (100,180.25) -- (100,190.25) (131,180.25) -- (131,190.25) (162,180.25) -- (162,190.25) (193,180.25) -- (193,190.25) (224,180.25) -- (224,190.25) (255,180.25) -- (255,190.25) (286,180.25) -- (286,190.25) (317,180.25) -- (317,190.25) (348,180.25) -- (348,190.25) (379,180.25) -- (379,190.25) (410,180.25) -- (410,190.25) (441,180.25) -- (441,190.25) (472,180.25) -- (472,190.25) (503,180.25) -- (503,190.25) (534,180.25) -- (534,190.25) (565,180.25) -- (565,190.25) (596,180.25) -- (596,190.25) (627,180.25) -- (627,190.25) (64,154.25) -- (74,154.25) (64,61.25) -- (74,61.25);

   \draw (100.07,154.33) circle[radius=2.75pt];
   \draw (130.40,154.33) circle[radius=2.75pt];
   \draw (193.07,154.33) circle[radius=2.75pt];
   \draw (223.90,154.33) circle[radius=2.75pt];
   \draw (255.23,154.33) circle[radius=2.75pt];
   \draw (316.07,154.33) circle[radius=2.75pt];
   \draw (347.73,154.33) circle[radius=2.75pt];
   \draw (410.07,154.33) circle[radius=2.75pt];
   \draw (471.73,154.33) circle[radius=2.75pt];
   \draw (502.57,154.33) circle[radius=2.75pt];
   \draw (564.23,154.33) circle[radius=2.75pt];
   \draw (595.07,154.33) circle[radius=2.75pt];
   \draw (627.07,154.33) circle[radius=2.75pt];
   
   \draw (100.07,61.33) circle[radius=2.75pt];
   \draw (130.40,61.33) circle[radius=2.75pt];
   \draw (193.07,61.33) circle[radius=2.75pt];
   \draw (223.90,61.33) circle[radius=2.75pt];
   \draw (255.23,61.33) circle[radius=2.75pt];
   \draw (316.07,61.33) circle[radius=2.75pt];
   \draw (347.73,61.33) circle[radius=2.75pt];
   \draw (410.07,61.33) circle[radius=2.75pt];
   \draw (471.73,61.33) circle[radius=2.75pt];
   \draw (502.57,61.33) circle[radius=2.75pt];
   \draw (564.23,61.33) circle[radius=2.75pt];
   \draw (595.07,61.33) circle[radius=2.75pt];
   \draw (627.07,61.33) circle[radius=2.75pt];

   \draw [->,>={Triangle[scale=.6]},line width=.35mm,dash pattern={on 0.84pt off 1.81pt}] (103.80,61.33) -- (125.67,61.33);
   \draw [->,>={Triangle[scale=.6]},line width=.35mm,dash pattern={on 0.84pt off 1.81pt}] (134.13,61.33) -- (188.33,61.33);
   \draw [->,>={Triangle[scale=.6]},line width=.35mm,dash pattern={on 0.84pt off 1.81pt}] (196.80,61.33) -- (219.17,61.33);
   \draw [->,>={Triangle[scale=.6]},line width=.35mm,dash pattern={on 0.84pt off 1.81pt}] (228.13,61.33) -- (250.50,61.33);
   \draw [->,>={Triangle[scale=.6]},line width=.35mm,dash pattern={on 0.84pt off 1.81pt}] (258.97,61.33) -- (311.33,61.33);
   \draw [->,>={Triangle[scale=.6]},line width=.35mm,dash pattern={on 0.84pt off 1.81pt}] (319.80,61.33) -- (342.17,61.33);
   \draw [->,>={Triangle[scale=.6]},line width=.35mm,dash pattern={on 0.84pt off 1.81pt}] (352.47,61.33) -- (405.67,61.33);
   \draw [->,>={Triangle[scale=.6]},line width=.35mm,dash pattern={on 0.84pt off 1.81pt}] (413.80,61.33) -- (467.00,61.33);
   \draw [->,>={Triangle[scale=.6]},line width=.35mm,dash pattern={on 0.84pt off 1.81pt}] (475.47,61.33) -- (497.83,61.33);
   \draw [->,>={Triangle[scale=.6]},line width=.35mm,dash pattern={on 0.84pt off 1.81pt}] (506.30,61.33) -- (559.50,61.33);
   \draw [->,>={Triangle[scale=.6]},line width=.35mm,dash pattern={on 0.84pt off 1.81pt}] (567.97,61.33) -- (590.33,61.33);
   \draw [->,>={Triangle[scale=.6]},line width=.35mm,dash pattern={on 0.84pt off 1.81pt}] (598.80,61.33) -- (622.17,61.33);

   \draw [->,>={Triangle[scale=.6]},line width=.35mm,dash pattern={on 0.84pt off 1.81pt}] (599.97,154.33) -- (622.33,154.33);
   \draw [->,>={Triangle[scale=.6]},line width=.35mm,dash pattern={on 0.84pt off 1.81pt}] (567.97,154.33) -- (590.33,154.33);
   \draw [->,>={Triangle[scale=.6]},line width=.35mm,dash pattern={on 0.84pt off 1.81pt}] (506.63,154.33) -- (559.50,154.33);
   \draw [->,>={Triangle[scale=.6]},line width=.35mm,dash pattern={on 0.84pt off 1.81pt}] (475.47,154.33) -- (497.83,154.33);
   \draw [->,>={Triangle[scale=.6]},line width=.35mm,dash pattern={on 0.84pt off 1.81pt}] (413.47,154.33) -- (467.00,154.33);
   \draw [->,>={Triangle[scale=.6]},line width=.35mm,dash pattern={on 0.84pt off 1.81pt}] (352.47,154.33) -- (405.33,154.33);
   \draw [->,>={Triangle[scale=.6]},line width=.35mm,dash pattern={on 0.84pt off 1.81pt}] (320.63,154.33) -- (343.00,154.33);
   \draw [->,>={Triangle[scale=.6]},line width=.35mm,dash pattern={on 0.84pt off 1.81pt}] (258.97,154.33) -- (311.33,154.33);
   \draw [->,>={Triangle[scale=.6]},line width=.35mm,dash pattern={on 0.84pt off 1.81pt}] (227.63,154.33) -- (250.00,154.33);
   \draw [->,>={Triangle[scale=.6]},line width=.35mm,dash pattern={on 0.84pt off 1.81pt}] (196.80,154.33) -- (219.17,154.33);
   \draw [->,>={Triangle[scale=.6]},line width=.35mm,dash pattern={on 0.84pt off 1.81pt}] (134.13,154.33) -- (188.33,154.33);
   \draw [->,>={Triangle[scale=.6]},line width=.35mm,dash pattern={on 0.84pt off 1.81pt}] (103.80,154.33) -- (126.17,154.33);

   \draw ( 62.60, 28.00) node [color={rgb, 255:red, 74; green, 144; blue, 226}, opacity=1, rotate=-0.49] {\large $\matheuler{S}$};
   \draw (674.60,186.33) node [color={rgb, 255:red, 74; green, 144; blue, 226}, opacity=1, rotate=-0.49] {\large $\matheuler{T}$};

\end{tikzpicture}
  }
  \caption{Example of an instance of the~\shortname modeled as a space-time network.~In this example, there are two stations, with nodes at the same station are grouped horizontally.~Connecting arcs are displayed as dotted arcs and demand arcs are displayed as full arcs.~All demand arcs with the same color stand for corresponding demands of the same customer.}
  \label{fig:network}
\end{figure}

When a demand is fulfilled by a vehicle from station~\(\outgoing{s}\) to~\(\incoming{s}\), starting at time~\(t_{i}\) and ending at time~\(t_{j}\), a unit flow is sent from node~\(\outgoing{s}_{i}\) to~\(\incoming{s}_{j}\).~Note that a unit flow along with a connecting arc~\mbox{\(e = (s_{i - 1}, s_{i})\)} corresponds to parking~(and possibly recharging the vehicle) at station~\(s\) from time~\(t_{i-1}\) until time~\(t_{i}\), in such a way that the recharging period means a left-closed and right-open interval.

Consider the following notation for sets of incoming and outgoing arcs for a node~\mbox{\(s_{i} \in \matheuler{N}_{s}\)}, \mbox{\(s \in \matheuler{S}\)}.~Let~\mbox{\(\incoming{\matheuler{E}}(s_{i}) = \{(s_{j}, s_{i}) \in \matheuler{E}_{c}, \forall c \in \matheuler{C}\}\)} and~\mbox{\(\outgoing{\matheuler{E}}(s_{i}) = \{(s_{i}, s_{j}) \in \matheuler{E}_{c}, \forall c \in \matheuler{C}\}\)} respectively be sets of incoming and outgoing demand arcs associated with a node~\(s_{i}\).~Additionally, we use the notation~\mbox{\(\incoming{e}(s_{i}) = (s_{i-1}, s_{i}) \in \matheuler{E}_{s}\)} and~\mbox{\(\outgoing{e}(s_{i}) = (s_{i}, s_{i + 1}) \in \matheuler{E}_{s}\)} to denote the incoming and outgoing connecting arcs associated with a node~\(s_{i}\), respectively.

\subsubsection*{Formulation \RemoveSpaces{\ref*{formulation:first}}}

For each~\mbox{\(v \in \matheuler{V}\)}, there are three sets of flow decision variables.~One of these sets is obtained by defining a binary variable~\(x^{v}_{e}\) for each demand arc~\(e\), which takes the value~\(1\) if and only if the demand associated with arc~\(e\) is fulfilled by vehicle~\(v\).~The other two sets are obtained by respectively defining binary variables~\(y^{v}_{e}\) and~\(z^{v}_{e}\) for each connecting arc~\(e\).~For each pair of variables~\(y^{v}_{e}\) and~\(z^{v}_{e}\), at most one of these variables can take the value~\(1\), if and only if the vehicle~\(v\) is parked at the station during the time period associated with the arc~\(e\).~Additionally,~\(y^{v}_{e}\) indicates that the vehicle~\(v\) is in a parking space without a charging facility, while~\(z^{v}_{e}\) indicates that the vehicle~\(v\) is in a parking space with a charging facility.

We define a set of continuous variables to represent the remaining energy of the battery over time, for each~\mbox{\(v \in \matheuler{V}\)}, obtained by defining a variable~\(\ell^{v}_{i}\), for all~\(i \in \matheuler{I}\).~We also define a binary variable~\(w_{c}\) indicating whether customer~\(c\) is served, which takes the value~\(1\) if and only if all demands of the customer~\(c\) are fulfilled, for each~\mbox{\(c \in \matheuler{C}\)}.

Using the notation above, the~\shortname is formulated as the following mixed-integer linear programming problem in a space-time network.

\noindent\RemoveSpaces{\makerefF{formulation:first}}
\begingroup
\allowdisplaybreaks
\small
\begin{alignat}{8}
\mbox{Maximize\,} & \sum\centerlap{c \in \matheuler{C}}\sum\rightlap{(%
                               \outgoing{s}\hspace*{-.3mm},%
                               t_{\hspace*{-.3mm}i}\hspace*{-.3mm},%
                               \incoming{s}\hspace*{-.3mm},%
                               t_{\hspace*{-.3mm}j}\!) \in \matheuler{D}_{c}\\} 
                \left( t_{j}-t_{i} \right) w_{c} & & & \makeref{EVSP1:ObjectiveFunction}\\[-7pt]
\mbox{Subject to}\! & \notag \\[-4pt]
            & \sum\centerlap{v \in \matheuler{V}} x^{v}_{e} = w_{c}
            & \forall c \in \matheuler{C},\,
              \forall e \in \matheuler{E}_{c}
            &
            & \makeref{EVSP1:DemandConstraints}\\
            & \sum\rightlap{e \in \incoming{\matheuler{E}}\spacetimezeroskips{s}{i}} x^{v}_{e} +
               y^{v}_{\incoming{e}\spacetimezeroskips{s}{i}} +
               z^{v}_{\incoming{e}\spacetimezeroskips{s}{i}} =
              \sum\rightlap{e \in \outgoing{\matheuler{E}}\spacetimezeroskips{s}{i}} x^{v}_{e} +
               y^{v}_{\outgoing{e}\spacetimezeroskips{s}{i}} + 
               z^{v}_{\outgoing{e}\spacetimezeroskips{s}{i}}
            & \forall v \in \matheuler{V},\,
              \forall s \in \matheuler{S},\,
              \forall {i} \in \matheuler{I} \!\setminus\! \{0, m\}
            & 
            & \makeref{EVSP1:FlowConservationEquation}\\
            & ~y^{v}_{\incoming{e}\spacetimezeroskips{s}{i}} \leq
               y^{v}_{\outgoing{e}\spacetimezeroskips{s}{i}} +
              \sum\rightlap{e \in \outgoing{\matheuler{E}}\spacetimezeroskips{s}{i}} x^{v}_{e}
            & \forall v \in \matheuler{V},\,
              \forall s \in \matheuler{S},\,
              \forall {i} \in \matheuler{I} \!\setminus\! \{0, m\}
            & 
            & \makeref{EVSP1:ConnectingSequence1}\\
            & ~z^{v}_{\incoming{e}\spacetimezeroskips{s}{i}} \leq
               z^{v}_{\outgoing{e}\spacetimezeroskips{s}{i}} +
              \sum\rightlap{e \in \outgoing{\matheuler{E}}\spacetimezeroskips{s}{i}} x^{v}_{e}
            & \forall v \in \matheuler{V},\,
              \forall s \in \matheuler{S},\,
              \forall {i} \in \matheuler{I} \!\setminus\! \{0, m\}
            &  
            & \makeref{EVSP1:ConnectingSequence2}\\
            & \sum\centerlap{v \in \matheuler{V}}\,
              \sum\rightlap{ e \in \incoming{\matheuler{E}}\spacetimezeroskips{s}{i}} x^{v}_{e} + 
              \sum\centerlap{v \in \matheuler{V}}\!
              \left( y^{v}_{\incoming{e}\spacetimezeroskips{s}{i}} \!+ 
               z^{v}_{\incoming{e}\spacetimezeroskips{s}{i}} \right) 
              \leq \textsf{C}_{s}
            & \forall s \in \matheuler{S},\,
              \forall i \in \matheuler{I} : \exists e \in \incoming{\matheuler{E}}(s_{i})
            &
            & \makeref{EVSP1:StationCapacity}\\
            & \sum\rightlap{v^{_\prime}\hspace*{-.2mm} \in \matheuler{V} 
                                                       \setminus\{\hspace*{-.2mm}v\hspace*{-.2mm}\}} 
               y^{v^{_\prime}}_{\incoming{e}\spacetimezeroskips{s}{i}} + 
               y^{v}_{\outgoing{e}\spacetimezeroskips{s}{i}} 
              \leq \textsf{C}_{s} - \textsf{R}_{s}
            & \hspace*{1.1cm}
              \forall s \in \matheuler{S},\,
              \forall v \in \matheuler{V},\,
              \forall i \in \matheuler{I} : \exists e \in \incoming{\matheuler{E}}(s_{i})
            & 
            & \makeref{EVSP1:ChargingPointCapacity1}\\
            & \sum\rightlap{v^{_\prime}\hspace*{-.2mm} \in \matheuler{V} 
                                                       \setminus\{\hspace*{-.2mm}v\hspace*{-.2mm}\}} 
               z^{v^{_\prime}}_{\incoming{e}\spacetimezeroskips{s}{i}} + 
               z^{v}_{\outgoing{e}\spacetimezeroskips{s}{i}} 
              \leq \textsf{R}_{s}
            & \forall s \in \matheuler{S},\,
              \forall v \in \matheuler{V},\,
              \forall i \in \matheuler{I} : \exists e \in \incoming{\matheuler{E}}(s_{i})
            & 
            & \makeref{EVSP1:ChargingPointCapacity2}\\
            & ~\ell^{v}_{i} \leq  \ell^{v}_{i-1} +
              \sum\centerlap{s \in \matheuler{S}}\!
              \raisebox{-2pt}{\bigg(}\!
                  \textsf{E}_{\incoming{e}\spacetimezeroskips{s}{i}} 
                   z^{v}_{\incoming{e}\spacetimezeroskips{s}{i}} -
                  \sum\rightlap{e \in \incoming{\matheuler{E}}\spacetimezeroskips{s}{i}}
                  \textsf{E}_{e} x^{v}_{e}\!
              \raisebox{-2pt}{\bigg)}
            & \forall v \in \matheuler{V},\,
              \forall i \in \matheuler{I} \!\setminus\! \{0\}
            &
            & \makeref{EVSP1:RechargingConstraints}\\
            & ~\ell^{v}_{i} \leq \textsf{L}
            & \forall v \in \matheuler{V},\,
              \forall i \in \matheuler{I} \!\setminus\! \{0\}
            &
            & \makeref{EVSP1:RechargingConstraints2}\\
            & ~\ell^{v}_{0} = \textsf{L}_{v}
            & \forall v \in \matheuler{V}
            & 
            & \makeref{EVSP1:InitialBatteryStateOfCharge}\\
            & ~y^{v}_{\outgoing{e}\spacetimezeroskips{s}{0}} = 1, 
              ~z^{v}_{\outgoing{e}\spacetimezeroskips{s}{0}} = 0
            & \forall s \in \matheuler{S},\,
              \forall v \in \matheuler{V}_{s} \!\!\setminus\! \matheuler{V}^{\prime}_{s}
            &
            & \makeref{EVSP1:InitialDistribution1}\\
            & ~y^{v}_{\outgoing{e}\spacetimezeroskips{s}{0}} = 0,
              ~z^{v}_{\outgoing{e}\spacetimezeroskips{s}{0}} = 1
            & \forall s \in \matheuler{S},\,
              \forall v \in \matheuler{V}^{\prime}_{s}
            &
            & \makeref{EVSP1:InitialDistribution2}\\
            & ~y^{v}_{\outgoing{e}\spacetimezeroskips{s}{0}} = 0,
              ~z^{v}_{\outgoing{e}\spacetimezeroskips{s}{0}} = 0
            & \forall s \in \matheuler{S},\,
              \forall v \in \matheuler{V} \!\setminus\!\! \matheuler{V}_{s}
            & 
            & \makeref{EVSP1:InitialDistribution3}\\
            & ~\ell^{v}_{i} \in \mathbb{R}_{\geq 0}
            & \forall v \in \matheuler{V},\, 
              \forall i \in \matheuler{I}
            &
            & \makeref{EVSP1:NonNegativeConstraints}\\
            & ~w_{c} \in \{0, 1\}
            & \forall c \in \matheuler{C}
            &
            & \makeref{EVSP1:IntegralityConstraints1}\\
            & ~x^{v}_{e} \in \{0, 1\}
            & \forall v \in \matheuler{V},\, 
              \forall e \in \bigcup\rightlap{c \in \matheuler{C}} \matheuler{E}_{c}
            &
            & \makeref{EVSP1:IntegralityConstraints2}\\
            & ~y^{v}_{e} \in \{0, 1\},\, z^{v}_{e} \in \{0, 1\}
            & \forall v \in \matheuler{V},\, 
              \forall e \in \bigcup\rightlap{s \in \matheuler{S}} \matheuler{E}_{s}
            &
            & \makeref{EVSP1:IntegralityConstraints4}
\end{alignat}
\endgroup

The objective function~\eqrefsmallsize{EVSP1:ObjectiveFunction} maximizes the sum of the rental times.~Constraint set~\eqrefsmallsize{EVSP1:DemandConstraints} makes sure that each demand for every served customer must be fulfilled by a single vehicle.~Constraint set~\eqrefsmallsize{EVSP1:FlowConservationEquation} is the flow conservation equations for the vehicles at each station and time instant, and constraint sets~\eqrefsmallsize{EVSP1:ConnectingSequence1} and~\eqrefsmallsize{EVSP1:ConnectingSequence2} impose that each vehicle remains parked~(and recharging) in the same parking space until it is used by the next customer.~Constraint set~\eqrefsmallsize{EVSP1:StationCapacity} ensures that the number of vehicles parked at each station at any time instant does not exceed the capacity.~It also ensures that each arriving vehicle can only be parked if there is an empty parking space beforehand.~Similarly, constraint sets~\eqrefsmallsize{EVSP1:ChargingPointCapacity1} and~\eqrefsmallsize{EVSP1:ChargingPointCapacity2} impose that parking spaces unequipped and equipped with charging facilities must be empty before they are occupied, respectively.~Note that only at the last time, although the capacity of the stations is guaranteed, there may be a surplus of vehicles parked in parking spaces equipped~(or unequipped) with charging facilities.~Since the vehicle does not need to be plugged into a charging facility after fulfilling the last demand, this surplus of vehicles does not compromise the feasibility of the solution.~However, to avoid both cases of surplus, the additional inequalities~\eqrefsmallsize{EVSP1:Inequality1} and~\eqrefsmallsize{EVSP1:Inequality2} can be considered in the formulation.
\begingroup
\small
\begin{alignat}{8}
   \hphantom{\mbox{Subject to}}\! 
   & \sum\centerlap{v \in \matheuler{V}} y^{v}_{\outgoing{e}\spacetimezeroskips{s}{m}} 
     \leq \textsf{C}_{s} - \textsf{R}_{s}
   & \hspace*{8.3cm}
     \forall s \in \matheuler{S}
   &
   & \makeref{EVSP1:Inequality1}\\
   & \sum\centerlap{v \in \matheuler{V}} z^{v}_{\outgoing{e}\spacetimezeroskips{s}{m}}
     \leq \textsf{R}_{s}
   & \forall s \in \matheuler{S}
   &
   & \makeref{EVSP1:Inequality2}
\end{alignat}
\endgroup

Constraint sets~\eqrefsmallsize{EVSP1:RechargingConstraints} and~\eqrefsmallsize{EVSP1:InitialBatteryStateOfCharge} refer to the energy stored in each vehicle's battery at any time.~More precisely, constraint set~\eqrefsmallsize{EVSP1:RechargingConstraints} captures the battery consumption and recharging whereas constraint set~\eqrefsmallsize{EVSP1:RechargingConstraints2} ensures that the electric energy recharged does not exceed the battery capacity, and constraint set~\eqrefsmallsize{EVSP1:InitialBatteryStateOfCharge} indicates the initial energy amount stored in the battery of each vehicle at the beginning of the planning horizon.~Constraint sets~\eqrefsmallsize{EVSP1:InitialDistribution1} to~\eqrefsmallsize{EVSP1:InitialDistribution3} indicate the initial state of the vehicles at each station.~Finally, expressions~\eqrefsmallsize{EVSP1:NonNegativeConstraints} to \eqrefsmallsize{EVSP1:IntegralityConstraints4} are the domain constraints for the variables.

An important practical advantage of this formulation is that the integrality of the binary variables~\(w_{c}\) and either~\(y^{v}_{e}\) or~\(z^{v}_{e}\) can be relaxed to be continuous between \(0\) and \(1\), since these variables are guaranteed to take a binary value by constraint sets~\eqrefsmallsize{EVSP1:DemandConstraints} and~\eqrefsmallsize{EVSP1:FlowConservationEquation}, respectively.~This can be proven by induction on the time instants for which constraint sets~\eqrefsmallsize{EVSP1:InitialDistribution1}, \eqrefsmallsize{EVSP1:InitialDistribution2} and~\eqrefsmallsize{EVSP1:InitialDistribution3} express the base cases.

\subsubsection*{Formulation \RemoveSpaces{\ref*{formulation:first_splitted}}}

A stronger formulation than~\RemoveSpaces{\ref{formulation:first}} can be derived by splitting each variable~\(x^{v}_{e}\) into four other variables, as follows.

For simplicity, we assume that the number of charging facilities~\(\textsf{R}_{s}\) is positive, but smaller than the capacity~\(\textsf{C}_{s}\), for each station~\mbox{\(s \in \matheuler{S}\)}.~There are four possibilities for each variable~\(x^{v}_{e}\).~These possibilities are based on outgoing from a parking space either unequipped or equipped with a charging facility and incoming to a parking space either unequipped or equipped with a charging facility.

Let~\mbox{\(\matheuler{K}= \{\textsf{UU}, \textsf{UE}, \textsf{EU}, \textsf{EE}\}\)} denote the set of possibilities for the split of the variable~\(x^{v}_{e}\).~To simplify the indexing, let~\mbox{\(\outgoing{\matheuler{Y}} = \{\textsf{UU}, \textsf{UE}\}\)} and~\mbox{\(\incoming{\matheuler{Y}} = \{\textsf{UU}, \textsf{EU}\}\)} be the sets containing all possibilities associated with outgoing from and incoming to a parking space unequipped with a charging facility, respectively.~Analogously, let~\mbox{\(\outgoing{\matheuler{Z}} = \{\textsf{EU}, \textsf{EE}\}\)} and~\mbox{\(\incoming{\matheuler{Z}} = \{\textsf{UE}, \textsf{EE}\}\)} be the sets containing all possibilities associated with outgoing from and incoming to a parking space equipped with a charging facility, respectively.

For each vehicle~\mbox{\(v \in \matheuler{V}\)} and each demand arc~\mbox{\(e \in \matheuler{E}_{c}\)}, for all~\mbox{\(c \in \matheuler{C}\)}, we define four binary variables~\mbox{\(x^{v}_{e\hspace*{-.2mm},\hspace*{-.2mm}k}\)}, for all~\mbox{\(k \in \matheuler{K}\)}, in which at most one of these variables can take the value~\(1\) if and only if the demand associated with arc~\(e\) is fulfilled by vehicle~\(v\).~In addition, each binary variable~\mbox{\(x^{v}_{e\hspace*{-.2mm},\hspace*{-.2mm}k}\)}, for all~\mbox{\(\smash{k \in \outgoing{\matheuler{Y}}}\)}, takes the value~\(1\) if and only if vehicle~\(v\) was in a parking space unequipped with a charging facility, while each binary variable~\mbox{\(x^{v}_{e\hspace*{-.2mm},\hspace*{-.2mm}k}\)}, for all~\mbox{\(k \in \outgoing{\matheuler{Z}}\)}, takes the value~\(1\) if and only if the vehicle~\(v\) was in a parking space equipped with a charging facility.~We use the remaining sets of variables included in the formulation~\RemoveSpaces{\ref{formulation:first}}, namely the variables~\(y^{v}_{e}\), \(z^{v}_{e}\), \(\ell^{v}_{i}\), and~\(w_{c}\). 


By introducing binary variables~\mbox{\(x^{v}_{e\hspace*{-.2mm},\hspace*{-.2mm}k}\)}, for~\mbox{\(k \in \matheuler{K}\)}, the~\shortname is reformulated as the following mixed-integer linear programming problem with split variables.

\noindent\RemoveSpaces{\makerefS{formulation:first_splitted}}
\begingroup
\allowdisplaybreaks
\small
\begin{alignat}{8}
\mbox{Maximize\,} & \sum\centerlap{c \in \matheuler{C}}\sum\rightlap{(%
                               \outgoing{s}\hspace*{-.3mm},%
                               t_{\hspace*{-.3mm}i}\hspace*{-.3mm},%
                               \incoming{s}\hspace*{-.3mm},%
                               t_{\hspace*{-.3mm}j}\!) \in \matheuler{D}_{c}\\} 
                \left( t_{j}-t_{i} \right) w_{c} & & & \makeref{EVSP1splitted:ObjectiveFunction}\\[-7pt]
\mbox{Subject to}\! & \notag \\[-4pt]
            & \sum\centerlap{v \in \matheuler{V}}\, 
              \sum\centerlap{k \in \matheuler{K}} x^{v}_{e\hspace*{-.2mm},\hspace*{-.2mm}k} = w_{c}
            & \forall c \in \matheuler{C},\,
              \forall e \in \matheuler{E}_{c}
            &
            & \makeref{EVSP1splitted:DemandConstraints}\\
            & \sum\centerlap{k \in \incoming{\matheuler{Y}}}\,
              \sum\rightlap{e \in \incoming{\matheuler{E}}\spacetimezeroskips{s}{i}}
               x^{v}_{e\hspace*{-.2mm},\hspace*{-.2mm}k} +
               y^{v}_{\incoming{e}\spacetimezeroskips{s}{i}} =
              \sum\centerlap{k \in \outgoing{\matheuler{Y}}}\,
              \sum\rightlap{e \in \outgoing{\matheuler{E}}\spacetimezeroskips{s}{i}}\,
               x^{v}_{e\hspace*{-.2mm},\hspace*{-.2mm}k} +
               y^{v}_{\outgoing{e}\spacetimezeroskips{s}{i}}
            & \forall v \in \matheuler{V},\,
              \forall s \in \matheuler{S},\,
              \forall {i} \in \matheuler{I} \!\setminus\! \{0, m\}
            &
            & \makeref{EVSP1splitted:FlowConservationEquation1}\\
            & \sum\centerlap{k \in \incoming{\matheuler{Z}}}\,
              \sum\rightlap{e \in \incoming{\matheuler{E}}\spacetimezeroskips{s}{i}}
               x^{v}_{e\hspace*{-.2mm},\hspace*{-.2mm}k} +
               z^{v}_{\incoming{e}\spacetimezeroskips{s}{i}} =
              \sum\centerlap{k \in \outgoing{\matheuler{Z}}}\,
              \sum\rightlap{e \in \outgoing{\matheuler{E}}\spacetimezeroskips{s}{i}}\,
               x^{v}_{e\hspace*{-.2mm},\hspace*{-.2mm}k} +
               z^{v}_{\outgoing{e}\spacetimezeroskips{s}{i}}
            & \forall v \in \matheuler{V},\,
              \forall s \in \matheuler{S},\,
              \forall {i} \in \matheuler{I} \!\setminus\! \{0, m\}
            & 
            & \makeref{EVSP1splitted:FlowConservationEquation2}\\
            & \sum\centerlap{v \in \matheuler{V}} \!\raisebox{-4pt}{\Bigg(}\!
                \sum\centerlap{k \in \incoming{\matheuler{Y}}}\,
                \sum\rightlap{ e \in \incoming{\matheuler{E}}\spacetimezeroskips{s}{i}}%
                 x^{v}_{e\hspace*{-.2mm},\hspace*{-.2mm}k} + 
                 y^{v}_{\incoming{e}\spacetimezeroskips{s}{i}}
              \!\raisebox{-4pt}{\Bigg)}
              \leq \textsf{C}_{s} - \textsf{R}_{s}
            & \forall s \in \matheuler{S},\,
              \forall i \in \matheuler{I} : \exists e \in \incoming{\matheuler{E}}(s_{i})
            &
            & \makeref{EVSP1splitted:StationCapacity}\\
            & \sum\centerlap{v \in \matheuler{V}} \!\raisebox{-4pt}{\Bigg(}\!
                \sum\centerlap{k \in \incoming{\matheuler{Z}}}\,
                \sum\rightlap{ e \in \incoming{\matheuler{E}}\spacetimezeroskips{s}{i}}%
                 x^{v}_{e\hspace*{-.2mm},\hspace*{-.2mm}k} + 
                 z^{v}_{\incoming{e}\spacetimezeroskips{s}{i}}
              \!\raisebox{-4pt}{\Bigg)}
              \leq \textsf{R}_{s}
            & \forall s \in \matheuler{S},\,
              \forall i \in \matheuler{I} : \exists e \in \incoming{\matheuler{E}}(s_{i})
            &
            & \makeref{EVSP1splitted:StationChargingFacilities}\\
            & ~\ell^{v}_{i} \leq \ell^{v}_{i-1} +
              \sum\centerlap{s \in \matheuler{S}}\!
              \raisebox{-4pt}{\Bigg(}\!
                  \textsf{E}_{\incoming{e}\spacetimezeroskips{s}{i}} 
                   z^{v}_{\incoming{e}\spacetimezeroskips{s}{i}} -
                  \sum\centerlap{k \in \matheuler{K}}\,
                  \sum\rightlap{e \in \incoming{\matheuler{E}}\spacetimezeroskips{s}{i}}
                  \textsf{E}_{e} x^{v}_{e\hspace*{-.2mm},\hspace*{-.2mm}k}\!
              \raisebox{-4pt}{\Bigg)}\!
              \hspace*{1.63cm}
            & \forall v \in \matheuler{V},\,
              \forall i \in \matheuler{I} \!\setminus\! \{0\}
            &
            & \makeref{EVSP1splitted:RechargingConstraints}\\
            & ~\ell^{v}_{i} \leq \textsf{L}
            & \forall v \in \matheuler{V},\,
              \forall i \in \matheuler{I} \!\setminus\! \{0\}
            &
            & \makeref{EVSP1splitted:RechargingConstraints2}
\end{alignat}
\endgroup

Note that, for brevity, we omit the constraint sets on the initial configuration of the system and the domain constraints on variables.

Constraint set~\eqrefsmallsize{EVSP1splitted:DemandConstraints} makes sure that each demand for every served customer must be fulfilled by a single vehicle.~Note that at most one of the four variables can be selected to fulfill the demand associated with the demand arc.~Constraint sets~\eqrefsmallsize{EVSP1splitted:FlowConservationEquation1} and~\eqrefsmallsize{EVSP1splitted:FlowConservationEquation2} are the flow conservation equations for the vehicles at each station and time instant.~Also, these constraint sets impose that each vehicle remains parked in the same parking space until it is used by the next customer.~Constraint sets~\eqrefsmallsize{EVSP1splitted:StationCapacity} and~\eqrefsmallsize{EVSP1splitted:StationChargingFacilities} ensure that the number of vehicles parked at each station at any time instant does not exceed the capacity, and that each arriving vehicle can only be parked if there is an empty parking space beforehand.~Constraint sets~\eqrefsmallsize{EVSP1splitted:RechargingConstraints} and~\eqrefsmallsize{EVSP1splitted:RechargingConstraints2} refer to the energy stored in each vehicle's battery along the planning horizon.

Not that the integrality of the binary variables~\(w_{c}\), \(z^{v}_{e}\) and~\(y^{v}_{e}\) can be relaxed to be continuous between \(0\) and \(1\), since these variables are guaranteed to take a binary value by the constraint sets~\eqrefsmallsize{EVSP1splitted:DemandConstraints}, \eqrefsmallsize{EVSP1splitted:FlowConservationEquation1} and~\eqrefsmallsize{EVSP1splitted:FlowConservationEquation2}, respectively.

\subsubsection*{Formulation strengths} 

We show that the linear relaxation of~\RemoveSpaces{\ref{formulation:first_splitted}} is at least as tight and sometimes tighter than the linear relaxation of~\RemoveSpaces{\ref{formulation:first}}.~For any such formulation~P, let~P$_{\textnormal{L}}$ denote its linear programming relaxation, and \mbox{let opt(P)} denote the value of its optimal objective value.

\begin{proposition}
  \label{proposition:EVSP1L_to_EVSP1-SL}
  Given any feasible solution of~\RemoveSpaces{\textnormal{\ref*{formulation:first_splitted}$_{\textnormal{L}}$}}, one can find a feasible solution of~\RemoveSpaces{\textnormal{\ref*{formulation:first}$_{\textnormal{L}}$}} with the same objective value.
\end{proposition}

{\noindent {\bf Proof.}~Consider a feasible solution for~\RemoveSpaces{\ref*{formulation:first_splitted}$_{\textnormal{L}}$}.~By making the sum of the values for~\mbox{\(x^{v}_{e\hspace*{-.2mm},\hspace*{-.2mm}k}\)}, for all~\mbox{\(k \in \matheuler{K}\)}, as the value of the variable~\(x^{v}_{e}\) in~\RemoveSpaces{\ref*{formulation:first}$_{\textnormal{L}}$}, and keeping the values of the remaining variables equal to the value of the same variable in the solution of~\RemoveSpaces{\ref*{formulation:first_splitted}$_{\textnormal{L}}$}~(\ie, direct mapping), we ensure that it has the same objective value.

Note that this mapping does not violate flow conservation constraints~\eqrefsmallsize{EVSP1:FlowConservationEquation}, because if we add constraints~\eqrefsmallsize{EVSP1splitted:FlowConservationEquation1} and~\eqrefsmallsize{EVSP1splitted:FlowConservationEquation2}, we obtain~\eqrefsmallsize{EVSP1:FlowConservationEquation}.~The additional parking constraints~\eqrefsmallsize{EVSP1:ConnectingSequence1} are satisfied since they can be derived from constraints~\eqrefsmallsize{EVSP1splitted:FlowConservationEquation1}, as follows:
\begingroup
\small
\begin{alignat*}{8}
   & \sum\centerlap{k \in \incoming{\matheuler{Y}}}\,
     \sum\rightlap{e \in \incoming{\matheuler{E}}\spacetimezeroskips{s}{i}}
       x^{v}_{e\hspace*{-.2mm},\hspace*{-.2mm}k} +
      y^{v}_{\incoming{e}\spacetimezeroskips{s}{i}} 
   & \,=\,
   & \sum\centerlap{k \in \outgoing{\matheuler{Y}}}\,
     \sum\rightlap{e \in \outgoing{\matheuler{E}}\spacetimezeroskips{s}{i}}\,
      x^{v}_{e\hspace*{-.2mm},\hspace*{-.2mm}k} +
     y^{v}_{\outgoing{e}\spacetimezeroskips{s}{i}}
   & \quad(\ref*{EVSP1splitted:FlowConservationEquation1})\\
   & 
   & \,\leq\,
   & \sum\centerlap{k \in \outgoing{\matheuler{Y}}}\,
     \sum\rightlap{e \in \outgoing{\matheuler{E}}\spacetimezeroskips{s}{i}}\,
       x^{v}_{e\hspace*{-.2mm},\hspace*{-.2mm}k} +
      y^{v}_{\outgoing{e}\spacetimezeroskips{s}{i}} +
     \sum\centerlap{k \in \outgoing{\matheuler{Z}}}\,
     \sum\rightlap{e \in \outgoing{\matheuler{E}}\spacetimezeroskips{s}{i}}\,
       x^{v}_{e\hspace*{-.2mm},\hspace*{-.2mm}k}
     \hspace*{2.4cm}
   & \\
   & 
   & \,=\,
   & \sum\centerlap{k \in \outgoing{\matheuler{K}}}\,
     \sum\rightlap{e \in \outgoing{\matheuler{E}}\spacetimezeroskips{s}{i}}\,
       x^{v}_{e\hspace*{-.2mm},\hspace*{-.2mm}k} +
      y^{v}_{\outgoing{e}\spacetimezeroskips{s}{i}}
   & \\
   & 
   & \,=\,
   & \sum\rightlap{e \in \outgoing{\matheuler{E}}\spacetimezeroskips{s}{i}}\,
       x^{v}_{e} +
      y^{v}_{\outgoing{e}\spacetimezeroskips{s}{i}}
   & \\[8pt]
   \hspace*{3cm}\Rightarrow~~
   & \sum\centerlap{k \in \incoming{\matheuler{Y}}}\,
     \sum\rightlap{e \in \incoming{\matheuler{E}}\spacetimezeroskips{s}{i}}\,
       x^{v}_{e\hspace*{-.2mm},\hspace*{-.2mm}k} +
      y^{v}_{\incoming{e}\spacetimezeroskips{s}{i}} 
   & \,\leq\,
   & \sum\rightlap{e \in \outgoing{\matheuler{E}}\spacetimezeroskips{s}{i}}\,
       x^{v}_{e} +
      y^{v}_{\outgoing{e}\spacetimezeroskips{s}{i}}
   & \\
   \Rightarrow~~
   & \hphantom{\sum\centerlap{k \in \incoming{\matheuler{Y}}}\,
     \sum\rightlap{e \in \incoming{\matheuler{E}}\spacetimezeroskips{s}{i}}\,
       x^{v}_{e\hspace*{-.2mm},\hspace*{-.2mm}k} +~~}
      y^{v}_{\incoming{e}\spacetimezeroskips{s}{i}} 
   & \,\leq\,
   & \sum\rightlap{e \in \outgoing{\matheuler{E}}\spacetimezeroskips{s}{i}}\,
       x^{v}_{e} +
      y^{v}_{\outgoing{e}\spacetimezeroskips{s}{i}}
   & \quad(\ref*{EVSP1:ConnectingSequence1})
\end{alignat*}
\endgroup

Analogously, constraints~\eqrefsmallsize{EVSP1:ConnectingSequence2} can be derived by~\eqrefsmallsize{EVSP1splitted:FlowConservationEquation2}.

The station capacity constraints~\eqrefsmallsize{EVSP1:StationCapacity} are satisfied since the addition of constraints~\eqrefsmallsize{EVSP1splitted:StationCapacity} and~\eqrefsmallsize{EVSP1splitted:StationChargingFacilities} gives constraint~\eqrefsmallsize{EVSP1:StationCapacity}, and the operational constraints~\eqrefsmallsize{EVSP1:ChargingPointCapacity1} are redundant with~\eqrefsmallsize{EVSP1splitted:StationCapacity}~if we algebraically manipulate the left-hand side, which can be rewritten as
\begingroup
\small
\begin{alignat*}{8}
   \hphantom{\mbox{Subject to}}\! 
   & \sum\centerlap{v^{_\prime}\! \in \matheuler{V}}\,
     \sum\centerlap{k \in \incoming{\matheuler{Y}}}\,
     \sum\rightlap{e \in \incoming{\matheuler{E}}\spacetimezeroskips{s}{i}}%
       x^{v^{_\prime}}_{e\hspace*{-.2mm},\hspace*{-.2mm}k} +
     \sum\rightlap{v^{_\prime}\! \in \matheuler{V}\setminus\{\hspace*{-.2mm}v\hspace*{-.2mm}\}}  
       y^{v^{_\prime}}_{\incoming{e}\spacetimezeroskips{s}{i}}
     + y^{v}_{\incoming{e}\spacetimezeroskips{s}{i}}
     \leq \textsf{C}_{s} - \textsf{R}_{s}
   & \hspace*{0.8cm}
     \hphantom{\forall} v \in \matheuler{V},\,
     \forall s \in \matheuler{S},\,
     \forall i \in \matheuler{I} : \exists e \in \incoming{\matheuler{E}}(s_{i})
   &
   & \hphantom{\quad(\ref*{EVSP1splitted:StationCapacity}^*)}
\end{alignat*}
\endgroup

\noindent by the fact that
\begingroup
\small
\begin{alignat*}{8}
   \hphantom{\mbox{Subject to}}\! 
   & \sum\centerlap{v^{_\prime}\! \in \matheuler{V}}\,
     \sum\centerlap{k \in \incoming{\matheuler{Y}}}\,
     \sum\rightlap{ e \in \incoming{\matheuler{E}}\spacetimezeroskips{s}{i}}%
       x^{v^{_\prime}}_{e\hspace*{-.2mm},\hspace*{-.2mm}k} +
     y^{v}_{\incoming{e}\spacetimezeroskips{s}{i}}
     \geq y^{v}_{\outgoing{e}\spacetimezeroskips{s}{i}}
   & \hspace*{2.75cm}
   \hphantom{\forall} v \in \matheuler{V},\,
     \forall s \in \matheuler{S},\,
     \forall i \in \matheuler{I} : \exists e \in \incoming{\matheuler{E}}(s_{i}),
   &
   & \hphantom{\quad(\ref*{EVSP1splitted:StationCapacity}^*)}
\end{alignat*}
\endgroup

\noindent which can be derived from~\eqrefsmallsize{EVSP1splitted:FlowConservationEquation1}.

In the same way, the operational constraints~\eqrefsmallsize{EVSP1:ChargingPointCapacity2} are redundant with~\eqrefsmallsize{EVSP1splitted:StationChargingFacilities} by applying algebraic manipulations on the left-hand side.~Clearly, all the remaining constraints are also satisfied.~\hfill \(\square\)
}



The following proposition shows that the optimal objective value of~\RemoveSpaces{\ref*{formulation:first_splitted}$_{\textnormal{L}}$} can be tighter than the value of~\RemoveSpaces{\ref*{formulation:first}$_{\textnormal{L}}$}. 

\begin{proposition}
  \label{proposition:EVSP1-SL_is_stronger}
  There are instances for which~\textnormal{opt(\RemoveSpaces{\ref*{formulation:first_splitted}$_{\textnormal{L}}$})} is smaller than~\textnormal{opt(\RemoveSpaces{\ref*{formulation:first}$_{\textnormal{L}}$})}.
\end{proposition}

{\noindent {\bf Proof.}~See the~\hyperref[appendix]{Appendix}.}~\hfill \(\square\)


\subsection{Heterogeneous space-time networks model}

As previously presented, after the construction of the network, a preprocessing phase can be performed to reduce the number of nodes and arcs, and hence the number of variables and constraints in the integer linear program.~First, two operations can be iteratively applied to the network, namely arc smoothing and node deleting.~These operations are applied to each node~\mbox{\(s_{i} \in \matheuler{N}_{s}\)}, \mbox{\(s \in \matheuler{S}\)}, that contains only the connecting arcs~\(\incoming{e}(s_{i})\) and~\(\outgoing{e}(s_{i})\), that is,~\(s_{i}\) has an in-degree and an out-degree of~\(1\).~An arc smoothing operation defined on~\mbox{\(\incoming{e}(s_{i}) = (s_{i-1}, s_{i})\)} and~\mbox{\(\outgoing{e}(s_{i}) = (s_{i}, s_{i+1})\)} is the process of removing the arcs~\(\incoming{e}(s_{i})\) and~\(\outgoing{e}(s_{i})\) while adding the arc~\mbox{\((s_{i-1}, s_{i+1})\)} to the network, and a node deleting operation defined on~\(s_{i}\) is the elementary process of removing the node from the network.

Another improvement can be obtained by fixing a certain number of variables.~This is based on the observation that for any vehicle~\mbox{\(v \in \matheuler{V}\)}, there may be demands for which we can ensure that they cannot be fulfilled by~\(v\).~In this sense, variables associated with these demands can be fixed to zero or removed from the model.

The following formulations are the result of combining of these two ideas.

\subsubsection*{Network representation and notation}

In order to reduce the number of variables, we consider a distinct space-time network for each vehicle.~To this end, we first determine the set of all demands that can be achievable by vehicle~\mbox{\(v \in \matheuler{V}\)}.~We say that a demand~\mbox{\(d = (\outgoing{s}, t_{i}, \incoming{s}, t_{j}) \in \matheuler{D}\)} is~\emph{achievable} by~\(v\) if, and only if, there is a feasible assignment plan for~\(v\) such that~\(d\) is fulfilled.

For each vehicle~\mbox{\(v \in \matheuler{V}\)}, let~\mbox{\(\matheuler{D}^{v}_{c} \subseteq \matheuler{D}_{c}\)} be the set of all demands of the customer~\mbox{\(c \in \matheuler{C}\)} that are achievable by~\(v\).~This set can be obtained by a label propagation algorithm for detecting achievable demands in the space-time network presented earlier.~The main idea behind this \mbox{algorithm} is that each node keeps the remaining energy for the vehicle~\(v\), which is propagated by visiting the nodes in increasing order of their time instants.~Initially, the remaining energy of the node where the vehicle~\(v\) is parked at time zero is~\(\textsf{L}_{v}\) and all other nodes have a remaining energy of~\(0\).~Then, at each iteration, the remaining energy of the current node is updated with the maximum energy to be obtained by a predecessor node, measured as the remaining energy of the predecessor node decreased by the required energy of demand, if the arc between the nodes is a demand arc; otherwise, increased by the amount of energy recharged at the station~(if there is any charging facility).~In addition, if the remaining energy decreased by the required energy of demand is non-negative, the demand arc is labeled as achievable.

A preprocessing step can be performed as follows after we have obtained all sets of achievable demands.~For each customer~\mbox{\(c \in \matheuler{C}\)} such that there is an unachievable demand by all vehicles in~\(\matheuler{V}\), remove all demands of~\(c\) from all sets of achievable demands.

Moreover, for each vehicle~\mbox{\(v \in \matheuler{V}\)}, let~\mbox{\(\matheuler{I}^{v}_{s} = \{ i \in \matheuler{I} : \exists (s, i) \in \Pi \} \cup \{ 0\}\)} be the set of indices for time instants at station~\mbox{\(s \in \matheuler{S}\)}, where~\(\Pi\) is the set of pairs of station and time index obtained as the union of all sets~\mbox{\(\{ (\outgoing{s}, i), (\incoming{s}, j) \}\)}, for each~\mbox{\((\outgoing{s}, t_{i}, \incoming{s}, t_{j}) \in \bigcup_{c \in \matheuler{C}} \matheuler{D}^{v}_{c}\)}; and let~\mbox{\(\matheuler{I}^{v} = \bigcup_{s \in \matheuler{S}} \matheuler{I}^{v}_{s}\)} denote the union of sets of indices for all stations.

In this model, the network for the vehicle~\(v\) is defined by a collection of node sets~\(\matheuler{N}^{v}_{s}\), for all~\mbox{\(s \in \matheuler{S}\)}, and a collection of arc sets~\(\matheuler{E}^{v}_{c}\) and~\(\matheuler{E}^{v}_{s}\), for all~\mbox{\(c \in \matheuler{C}\)} and~\mbox{\(s \in \matheuler{S}\)}, as follows: in each~\(\matheuler{N}^{v}_{s}\), there is a space-time node~\(s_{i}\), for all~\mbox{\(i \in \matheuler{I}^{v}_{s} \cup \{m\}\)}; in each \(\matheuler{E}^{v}_{c}\), there is a demand arc from node~\(\outgoing{s}_{i}\) to~\(\incoming{s}_{j}\), for all~\mbox{\((\outgoing{s}, t_{i}, \incoming{s}, t_{j}) \in \matheuler{D}^{v}_{c}\)}; and in each~\(\matheuler{E}^{v}_{s}\), there is a connecting arc from node~\(s_{i^{\prime}}\) to~\(s_{i}\), for all~\mbox{\(i \in \matheuler{I}^{v}_{s} \setminus\! \{0\} \cup  \{m\}\)} and~\mbox{\(i^{\prime} = \max\{ i^{\prime} \in \matheuler{I}^{v}_{s} : i^{\prime} < i \}\)}.

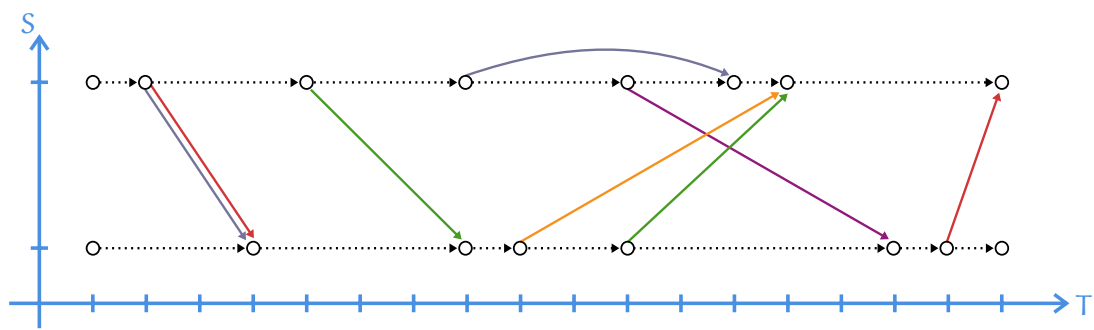
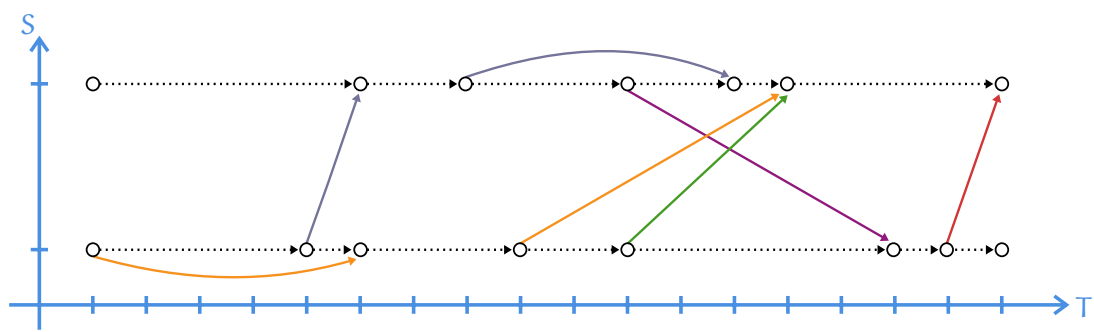
\begin{figure}[H]
  \centering
  \begin{subfigure}[t]{\textwidth}
    \resizebox{\textwidth}{!}{%
      \tikzset{every picture/.style={line width=0.75pt}}

\begin{tikzpicture}[x=0.75pt, y=0.75pt, yscale=-1, xscale=0.96]
   
   \draw [->,>={Triangle[scale=.6]},line width=.35mm,color={rgb, 255:red, 209; green, 53; blue, 53}, draw opacity=1] (133.5,62.63) -- (193.38,148.75);
   \draw [->,>={Triangle[scale=.6]},line width=.35mm,color=CadetBlue, draw opacity=1] (130.4,65.42) -- (188.81,149.89);
   \draw [->,>={Triangle[scale=.6]},line width=.35mm,color={rgb, 255:red, 70; green, 155; blue, 36}, draw opacity=1] (226.4,65.42) -- (313.89,149.5);
   \draw [->,>={Triangle[scale=.6]},line width=.35mm,color={rgb, 255:red, 145; green, 25; blue, 123}, draw opacity=1] (410.07,65.) -- (561.6,149.29);
   \draw [->,>={Triangle[scale=.6]},line width=.35mm,color=BurntOrange, draw opacity=1] (347.73,150.9) -- (498.1,66.89);
   \draw [->,>={Triangle[scale=.6]},line width=.35mm,color={rgb, 255:red, 70; green, 155; blue, 36}, draw opacity=1] (410.07,150.9) -- (502.91,67.57);
   \draw [->,>={Triangle[scale=.6]},line width=.35mm,color={rgb, 255:red, 209; green, 53; blue, 53}, draw opacity=1] (595.07,151.07) -- (625.52,67.06);
   \draw [->,>={Triangle[scale=.6]},line width=.35mm,color=CadetBlue, draw opacity=1] (316.07,57.57) .. controls (376.58,37.55) and (422.12,39.2) .. (469.09,57.23);

   \draw [color={rgb, 255:red, 74; green, 144; blue, 226}, draw opacity=1][line width=1.5] (51.5,185.25) -- (664.5,185.25) (69,36) -- (69,199.2) (657.5,180.25) -- (664.5,185.25) -- (657.5,190.25) (64,43) -- (69,36) -- (74,43) (100,180.25) -- (100,190.25) (131,180.25) -- (131,190.25) (162,180.25) -- (162,190.25) (193,180.25) -- (193,190.25) (224,180.25) -- (224,190.25) (255,180.25) -- (255,190.25) (286,180.25) -- (286,190.25) (317,180.25) -- (317,190.25) (348,180.25) -- (348,190.25) (379,180.25) -- (379,190.25) (410,180.25) -- (410,190.25) (441,180.25) -- (441,190.25) (472,180.25) -- (472,190.25) (503,180.25) -- (503,190.25) (534,180.25) -- (534,190.25) (565,180.25) -- (565,190.25) (596,180.25) -- (596,190.25) (627,180.25) -- (627,190.25) (64,154.25) -- (74,154.25) (64,61.25) -- (74,61.25);

   \draw (100.07,154.33) circle[radius=2.75pt];
   \draw (193.07,154.33) circle[radius=2.75pt];
   \draw (316.07,154.33) circle[radius=2.75pt];
   \draw (347.73,154.33) circle[radius=2.75pt];
   \draw (410.07,154.33) circle[radius=2.75pt];
   \draw (564.23,154.33) circle[radius=2.75pt];
   \draw (595.07,154.33) circle[radius=2.75pt];
   \draw (627.07,154.33) circle[radius=2.75pt];
   
   \draw (100.07,61.33) circle[radius=2.75pt];
   \draw (130.40,61.33) circle[radius=2.75pt];
   \draw (223.90,61.33) circle[radius=2.75pt];
   \draw (316.07,61.33) circle[radius=2.75pt];
   \draw (410.07,61.33) circle[radius=2.75pt];
   \draw (471.73,61.33) circle[radius=2.75pt];
   \draw (502.57,61.33) circle[radius=2.75pt];
   \draw (627.07,61.33) circle[radius=2.75pt];

   \draw [->,>={Triangle[scale=.6]},line width=.35mm,dash pattern={on 0.84pt off 1.81pt}] (103.80,61.33) -- (125.67,61.33);
   \draw [->,>={Triangle[scale=.6]},line width=.35mm,dash pattern={on 0.84pt off 1.81pt}] (134.13,61.33) -- (219.17,61.33);
   \draw [->,>={Triangle[scale=.6]},line width=.35mm,dash pattern={on 0.84pt off 1.81pt}] (228.13,61.33) -- (311.33,61.33);
   \draw [->,>={Triangle[scale=.6]},line width=.35mm,dash pattern={on 0.84pt off 1.81pt}] (319.80,61.33) -- (405.67,61.33);
   \draw [->,>={Triangle[scale=.6]},line width=.35mm,dash pattern={on 0.84pt off 1.81pt}] (413.80,61.33) -- (467.00,61.33);
   \draw [->,>={Triangle[scale=.6]},line width=.35mm,dash pattern={on 0.84pt off 1.81pt}] (475.47,61.33) -- (497.83,61.33);
   \draw [->,>={Triangle[scale=.6]},line width=.35mm,dash pattern={on 0.84pt off 1.81pt}] (506.30,61.33) -- (622.17,61.33);

   \draw [->,>={Triangle[scale=.6]},line width=.35mm,dash pattern={on 0.84pt off 1.81pt}] (599.97,154.33) -- (622.33,154.33);
   \draw [->,>={Triangle[scale=.6]},line width=.35mm,dash pattern={on 0.84pt off 1.81pt}] (567.97,154.33) -- (590.33,154.33);
   \draw [->,>={Triangle[scale=.6]},line width=.35mm,dash pattern={on 0.84pt off 1.81pt}] (413.47,154.33) -- (559.50,154.33);
   \draw [->,>={Triangle[scale=.6]},line width=.35mm,dash pattern={on 0.84pt off 1.81pt}] (352.47,154.33) -- (405.33,154.33);
   \draw [->,>={Triangle[scale=.6]},line width=.35mm,dash pattern={on 0.84pt off 1.81pt}] (320.63,154.33) -- (343.00,154.33);
   \draw [->,>={Triangle[scale=.6]},line width=.35mm,dash pattern={on 0.84pt off 1.81pt}] (196.80,154.33) -- (311.33,154.33);
   \draw [->,>={Triangle[scale=.6]},line width=.35mm,dash pattern={on 0.84pt off 1.81pt}] (103.80,154.33) -- (188.33,154.33);

   \draw ( 62.60, 28.00) node [color={rgb, 255:red, 74; green, 144; blue, 226}, opacity=1, rotate=-0.49] {\large $\matheuler{S}$};
   \draw (674.60,186.33) node [color={rgb, 255:red, 74; green, 144; blue, 226}, opacity=1, rotate=-0.49] {\large $\matheuler{T}$};

\end{tikzpicture}
    }
    \caption{An example for vehicles initially located at the top station.}
  \end{subfigure}
  
  \medskip
  \begin{subfigure}[t]{\textwidth}
    \resizebox{\textwidth}{!}{%
      \tikzset{every picture/.style={line width=0.75pt}}

\begin{tikzpicture}[x=0.75pt, y=0.75pt, yscale=-1, xscale=0.96]

   \draw [->,>={Triangle[scale=.6]},line width=.35mm,color=CadetBlue, draw opacity=1] (223.9,150.73) -- (236.56,116.47) -- (254.19,66.73);
   \draw [->,>={Triangle[scale=.6]},line width=.35mm,color=BurntOrange, draw opacity=1] (100.07,158.) .. controls (154.18,172.92) and (202.3,173.42) .. (252.92,159.56);
   \draw [->,>={Triangle[scale=.6]},line width=.35mm,color={rgb, 255:red, 145; green, 25; blue, 123}, draw opacity=1] (410.07,65.) -- (561.6,149.29);
   \draw [->,>={Triangle[scale=.6]},line width=.35mm,color=BurntOrange, draw opacity=1] (347.73,150.9) -- (498.1,66.89);
   \draw [->,>={Triangle[scale=.6]},line width=.35mm,color={rgb, 255:red, 70; green, 155; blue, 36}, draw opacity=1] (410.07,150.9) -- (502.91,67.57);
   \draw [->,>={Triangle[scale=.6]},line width=.35mm,color={rgb, 255:red, 209; green, 53; blue, 53}, draw opacity=1] (595.07,151.07) -- (625.52,67.06);
   \draw [->,>={Triangle[scale=.6]},line width=.35mm,color=CadetBlue, draw opacity=1] (316.07,57.57) .. controls (376.58,37.55) and (422.12,39.2) .. (469.09,57.23);
   
   \draw [color={rgb, 255:red, 74; green, 144; blue, 226}, draw opacity=1][line width=1.5] (51.5,185.25) -- (664.5,185.25) (69,36) -- (69,199.2) (657.5,180.25) -- (664.5,185.25) -- (657.5,190.25) (64,43) -- (69,36) -- (74,43) (100,180.25) -- (100,190.25) (131,180.25) -- (131,190.25) (162,180.25) -- (162,190.25) (193,180.25) -- (193,190.25) (224,180.25) -- (224,190.25) (255,180.25) -- (255,190.25) (286,180.25) -- (286,190.25) (317,180.25) -- (317,190.25) (348,180.25) -- (348,190.25) (379,180.25) -- (379,190.25) (410,180.25) -- (410,190.25) (441,180.25) -- (441,190.25) (472,180.25) -- (472,190.25) (503,180.25) -- (503,190.25) (534,180.25) -- (534,190.25) (565,180.25) -- (565,190.25) (596,180.25) -- (596,190.25) (627,180.25) -- (627,190.25) (64,154.25) -- (74,154.25) (64,61.25) -- (74,61.25);

   \draw (100.07,154.33) circle[radius=2.75pt];
   \draw (223.90,154.33) circle[radius=2.75pt];
   \draw (255.23,154.33) circle[radius=2.75pt];
   \draw (347.73,154.33) circle[radius=2.75pt];
   \draw (410.07,154.33) circle[radius=2.75pt];
   \draw (564.23,154.33) circle[radius=2.75pt];
   \draw (595.07,154.33) circle[radius=2.75pt];
   \draw (627.07,154.33) circle[radius=2.75pt];
   
   \draw (100.07,61.33) circle[radius=2.75pt];
   \draw (255.23,61.33) circle[radius=2.75pt];
   \draw (316.07,61.33) circle[radius=2.75pt];
   \draw (410.07,61.33) circle[radius=2.75pt];
   \draw (471.73,61.33) circle[radius=2.75pt];
   \draw (502.57,61.33) circle[radius=2.75pt];
   \draw (627.07,61.33) circle[radius=2.75pt];
   
   \draw [->,>={Triangle[scale=.6]},line width=.35mm,dash pattern={on 0.84pt off 1.81pt}] (103.80,61.33) -- (250.50,61.33);
   \draw [->,>={Triangle[scale=.6]},line width=.35mm,dash pattern={on 0.84pt off 1.81pt}] (258.97,61.33) -- (311.33,61.33);
   \draw [->,>={Triangle[scale=.6]},line width=.35mm,dash pattern={on 0.84pt off 1.81pt}] (319.80,61.33) -- (405.67,61.33);
   \draw [->,>={Triangle[scale=.6]},line width=.35mm,dash pattern={on 0.84pt off 1.81pt}] (413.80,61.33) -- (467.00,61.33);
   \draw [->,>={Triangle[scale=.6]},line width=.35mm,dash pattern={on 0.84pt off 1.81pt}] (475.47,61.33) -- (497.83,61.33);
   \draw [->,>={Triangle[scale=.6]},line width=.35mm,dash pattern={on 0.84pt off 1.81pt}] (506.30,61.33) -- (622.17,61.33);
   
   \draw [->,>={Triangle[scale=.6]},line width=.35mm,dash pattern={on 0.84pt off 1.81pt}] (599.97,154.33) -- (622.33,154.33);
   \draw [->,>={Triangle[scale=.6]},line width=.35mm,dash pattern={on 0.84pt off 1.81pt}] (567.97,154.33) -- (590.33,154.33);
   \draw [->,>={Triangle[scale=.6]},line width=.35mm,dash pattern={on 0.84pt off 1.81pt}] (413.47,154.33) -- (559.50,154.33);
   \draw [->,>={Triangle[scale=.6]},line width=.35mm,dash pattern={on 0.84pt off 1.81pt}] (352.47,154.33) -- (405.33,154.33);
   \draw [->,>={Triangle[scale=.6]},line width=.35mm,dash pattern={on 0.84pt off 1.81pt}] (258.97,154.33) -- (343.00,154.33);
   \draw [->,>={Triangle[scale=.6]},line width=.35mm,dash pattern={on 0.84pt off 1.81pt}] (227.63,154.33) -- (250.00,154.33);
   \draw [->,>={Triangle[scale=.6]},line width=.35mm,dash pattern={on 0.84pt off 1.81pt}] (103.80,154.33) -- (219.17,154.33);

   \draw ( 62.60, 28.00) node [color={rgb, 255:red, 74; green, 144; blue, 226}, opacity=1, rotate=-0.49] {\large $\matheuler{S}$};
   \draw (674.60,186.33) node [color={rgb, 255:red, 74; green, 144; blue, 226}, opacity=1, rotate=-0.49] {\large $\matheuler{T}$};

\end{tikzpicture}
    }
    \caption{An example for vehicles initially located at the bottom station.}
  \end{subfigure}
  \caption{Illustration of space-time networks constructed using the set of achievable \mbox{demands}, where each network represents vehicles located at one of two stations.}
  \label{fig:networks}
\end{figure}

To illustrate the network constructed as described above,~\Cref{fig:networks} gives two examples in which vehicles are initially at different stations~(the original network was shown in~\Cref{fig:network}).~Since the vehicles are identical, all vehicles initially located at the same station have the same set of achievable demands. 

We use the same notation as before for the amount~\(\textsf{E}_{e}\) of energy consumed by each demand arc~\mbox{\(e \in \matheuler{E}^{v}_{c}\)}, for all~\mbox{\(c \in \matheuler{C}\)}, and the amount~\(\textsf{E}_{e}\) of energy recharged at station~\(s\), in each connecting arc~\mbox{\(e \in \matheuler{E}^{v}_{s}\)}, for all~\mbox{\(s \in \matheuler{S}\)}.~Moreover, we use the same notation for the sets of incoming and outgoing demand arcs, denoted by~\(\incoming{\matheuler{E}}(s_{i})\) and~\(\outgoing{\matheuler{E}}(s_{i})\), respectively, and for the incoming and outgoing connecting arcs, denoted by~\(\incoming{e}(s_{i})\) and~\(\outgoing{e}(s_{i})\), from a node~\mbox{\(s_{i} \in \matheuler{N}^{v}_{s}\)}, for all~\mbox{\(s \in \matheuler{S}\)}.

Lastly, let~\mbox{\(\incoming{\matheuler{I}}_{s} = \bigcup_{v \in \matheuler{V}} \{ i \in \matheuler{I}^{v}_{s} : \exists e \in \incoming{\matheuler{E}}(s_{i}) \}\)} denote the set of indices for all possible arriving vehicles at station~\(s\), for each~\mbox{\(s \in \matheuler{S}\)}.

\subsubsection*{Formulation \RemoveSpaces{\ref*{formulation:second}}}

The decision variables used in this third formulation~(referred to as~\RemoveSpaces{\ref{formulation:second}}) have the same meaning as in the first formulation.~We start with the set of variables obtained by defining binary variables~\(x^{v}_{e}\) for all demand arc~\mbox{\(e \in \matheuler{E}^{v}_{c}\)}, for each~\mbox{\(v \in \matheuler{V}\)} and \mbox{\(c \in \matheuler{C}\)}, such that each~\(x^{v}_{e}\) takes the value~\(1\) if and only if the demand associated with arc~\(e\) is fulfilled by vehicle~\(v\).~Moreover, there are sets of binary variables~\(y^{v}_{e}\) and~\(z^{v}_{e}\) for all connecting arc~\mbox{\(e \in \matheuler{E}^{v}_{s}\)}, for each~\mbox{\(v \in \matheuler{V}\)} and~\mbox{\(s \in \matheuler{S}\)}, such that, for each pair, at most one of these variables can take the value~\(1\), if and only if the vehicle~\(v\) is parked at the station during the time period associated with~arc~\(e\).~Remembering that~\(z^{v}_{e}\) also indicates that the vehicle~\(v\) is parked within a parking space equipped with a charging facility.~To indicate the remaining energy of the battery over time, we define a set of continuous variables~\(\ell^{v}_{i}\) for each~\mbox{\(v \in \matheuler{V}\)}, \mbox{\(i \in \matheuler{I}^{v}\)}; and a binary variable~\(w_{c}\) indicates whether the customer~\(c\) is served, for each~\mbox{\(c \in \matheuler{C}\)}.~It takes the value~\(1\) if and only if all demands of customer~\(c\) are fulfilled.

The~\shortname is reformulated as the following mixed-integer linear programming problem by projecting the same set of constraints of the formulation~\RemoveSpaces{\ref{formulation:first}} over the heterogeneous space-time networks.~For clarity, this formulation is given next. 

\noindent\RemoveSpaces{\makerefF{formulation:second}}
\begingroup
\allowdisplaybreaks
\small
\begin{alignat}{8}
\mbox{Maximize\,} & \sum\centerlap{c \in \matheuler{C}}\sum\rightlap{(%
                               \outgoing{s}\hspace*{-.3mm},%
                               t_{\hspace*{-.3mm}i}\hspace*{-.3mm},%
                               \incoming{s}\hspace*{-.3mm},%
                               t_{\hspace*{-.3mm}j}\!) \in \matheuler{D}_{c}\\} 
                \left( t_{j}-t_{i} \right) w_{c} & & & \makeref{EVSP2:ObjectiveFunction}\\[-7pt]
\mbox{Subject to}\! & \notag \\[-4pt]
            & \sum\rightlap{v \in \matheuler{V}:
                \exists e {\scriptscriptstyle =}
                \doublespacetimezeroskips{\outgoing{s}}{i}{\incoming{s}}{j}\in\matheuler{E}^{v}_{c}
              }
               x^{v}_{e} = w_{c}
            & \forall c \in \matheuler{C},\,
              \forall (\outgoing{s}, t_{i}, \incoming{s}, t_{j}) \in \matheuler{D}_{c}
            &
            & \makeref{EVSP2:DemandConstraints}\\
            & \sum\rightlap{e \in \incoming{\matheuler{E}}\spacetimezeroskips{s}{i}} x^{v}_{e} +
              y^{v}_{\incoming{e}\spacetimezeroskips{s}{i}} + z^{v}_{\incoming{e}\spacetimezeroskips{s}{i}} =
              \sum\rightlap{e \in \outgoing{\matheuler{E}}\spacetimezeroskips{s}{i}} x^{v}_{e} +
              y^{v}_{\outgoing{e}\spacetimezeroskips{s}{i}} + z^{v}_{\outgoing{e}\spacetimezeroskips{s}{i}}
            & \forall v \in \matheuler{V},\,
              \forall s \in \matheuler{S},\,
              \forall i \in \matheuler{I}^{v} \!\setminus\!\! \{ 0, m \}
            & 
            & \makeref{EVSP2:FlowConservationEquation}\\
            & ~y^{v}_{\incoming{e}\spacetimezeroskips{s}{i}} \leq 
              y^{v}_{\outgoing{e}\spacetimezeroskips{s}{i}} +
              \sum\rightlap{e \in \outgoing{\matheuler{E}}\spacetimezeroskips{s}{i}} x^{v}_{e} 
            & \forall v \in \matheuler{V},\,
              \forall s \in \matheuler{S},\,
              \forall i \in \matheuler{I}^{v} \!\setminus\!\! \{ 0, m \}
            &
            & \makeref{EVSP2:ConnectingSequence1}\\
            & ~z^{v}_{\incoming{e}\spacetimezeroskips{s}{i}} \leq 
              z^{v}_{\outgoing{e}\spacetimezeroskips{s}{i}} +
              \sum\rightlap{e \in \outgoing{\matheuler{E}}\spacetimezeroskips{s}{i}} x^{v}_{e} 
            & \forall v \in \matheuler{V},\,
              \forall s \in \matheuler{S},\,
              \forall i \in \matheuler{I}^{v} \!\setminus\!\! \{ 0, m \}
            &
            & \makeref{EVSP2:ConnectingSequence2}\\
            & \sum\centerlap{v \in \matheuler{V}}\,
              \sum\rightlap{e \in \incoming{\matheuler{E}}\spacetimezeroskips{s}{i}} x^{v}_{e} +   
              \sum\rightlap{v \in \matheuler{V}:\exists\doublespacetimezeroskips{s}{j}{s}{k} \in 
                \matheuler{E}^{v}_{s}\!,\, j < i \leq k 
              }
              \left( y^{v}_{\incoming{e}\spacetimezeroskips{s}{k}} 
              + z^{v}_{\incoming{e}\spacetimezeroskips{s}{k}} \right) \leq \textsf{C}_{s}
            & \forall s \in \matheuler{S},\,
              \forall i \in \incoming{\matheuler{I}}_{s}
            &
            & \makeref{EVSP2:StationCapacity}\\
            & \sum\rightlap{v^{_\prime}\hspace*{-.2mm}
                \in \matheuler{V} \setminus\{\hspace*{-.2mm}v\hspace*{-.2mm}\}:
                \exists \doublespacetimezeroskips{s}{j}{s}{k}
                \in \matheuler{E}^{v^{_\prime}}_{s}\!\!,\, j < i \leq k 
              }
              y^{v^{_\prime}}_{\incoming{e}\spacetimezeroskips{s}{k}} +
              y^{v}_{\outgoing{e}\spacetimezeroskips{s}{j}} 
              \leq \textsf{C}_{s} - \textsf{R}_{s}
            & \hspace*{-2cm}
              \forall s \in \matheuler{S},\, 
              \forall i \in \incoming{\matheuler{I}}_{s},\,
              \forall v \in \matheuler{V} : 
              \exists (s_{j}\!, s_{k}) \in \matheuler{E}^{v}_{s}, j \leq i < k
            &
            & \makeref{EVSP2:ChargingPointCapacity1}\\
            & \sum\rightlap{v^{_\prime}\hspace*{-.2mm} 
                \in \matheuler{V} \setminus\{\hspace*{-.2mm}v\hspace*{-.2mm}\}:
                \exists \doublespacetimezeroskips{s}{j}{s}{k} 
                \in \matheuler{E}^{v^{_\prime}}_{s}\!\!,\, j < i \leq k 
              }
              z^{v^{_\prime}}_{\incoming{e}\spacetimezeroskips{s}{k}} 
              + z^{v}_{\outgoing{e}\spacetimezeroskips{s}{j}} \leq \textsf{R}_{s} 
            & \hspace*{-2cm}
              \forall s \in \matheuler{S},\, 
              \forall i \in \incoming{\matheuler{I}}_{s},\,
              \forall v \in \matheuler{V} : 
              \exists (s_{j}\!, s_{k}) \in \matheuler{E}^{v}_{s}, j \leq i < k
            &
            & \makeref{EVSP2:ChargingPointCapacity2}\\
            & ~\ell^{v}_{i} \leq \ell^{v}_{i^{\prime}} +
              \sum\centerlap{s \in \matheuler{S}}\!
              \raisebox{-2pt}{\bigg(}\!
                  \textsf{E}_{\incoming{e}\spacetimezeroskips{s}{i}} 
                  z^{v}_{\incoming{e}\spacetimezeroskips{s}{i}} -
                  \sum\rightlap{e \in \incoming{\matheuler{E}}\spacetimezeroskips{s}{i}}
                  \textsf{E}_{e} x^{v}_{e}\!
              \raisebox{-2pt}{\bigg)}
            & \hspace*{-.4cm}
              \forall v \in \matheuler{V},\,
              \forall i \in \matheuler{I}^{v} \!\setminus\!\! \{0\},\,
              i^{\prime} \!= \max \{ i^{\prime}\! \in \matheuler{I}^{v}\!: i^{\prime}\! < i\}
            &
            & \makeref{EVSP2:RechargingConstraints}\\
            & ~\ell^{v}_{i} \leq \textsf{L}
            & \forall v \in \matheuler{V},\,
              \forall i \in \matheuler{I}^{v} \!\setminus\!\! \{0\}
            &
            & \makeref{EVSP2:RechargingConstraints2}\\
            & ~\ell^{v}_{0} = \textsf{L}_{v}
            & \forall v \in \matheuler{V}
            &
            & \makeref{EVSP2:InitialBatteryStateOfCharge}\\
            & ~y^{v}_{\outgoing{e}\spacetimezeroskips{s}{0}} = 1, 
              ~z^{v}_{\outgoing{e}\spacetimezeroskips{s}{0}} = 0
            & \forall s \in \matheuler{S},\,
              \forall v \in \matheuler{V}_{s} \!\!\setminus\! \matheuler{V}^{\prime}_{s}
            &
            & \makeref{EVSP2:InitialDistribution1}\\
            & ~y^{v}_{\outgoing{e}\spacetimezeroskips{s}{0}} = 0,
              ~z^{v}_{\outgoing{e}\spacetimezeroskips{s}{0}} = 1
            & \forall s \in \matheuler{S},\,
              \forall v \in \matheuler{V}^{\prime}_{s}
            &
            & \makeref{EVSP2:InitialDistribution2}\\
            & ~y^{v}_{\outgoing{e}\spacetimezeroskips{s}{0}} = 0,
              ~z^{v}_{\outgoing{e}\spacetimezeroskips{s}{0}} = 0
            & \forall s \in \matheuler{S},\,
              \forall v \in \matheuler{V} \!\setminus\!\! \matheuler{V}_{s}
            & 
            & \makeref{EVSP2:InitialDistribution3}\\
            & ~\ell^{v}_{i} \in \mathbb{R}_{\geq 0}
            & \forall v \in \matheuler{V},\, 
              \forall i \in \matheuler{I}^{v}
            &
            & \makeref{EVSP2:NonNegativeConstraints}\\
            & ~w_{c} \in \{0, 1\}
            & \forall c \in \matheuler{C}
            &
            & \makeref{EVSP2:IntegralityConstraints1}\\
            & ~x^{v}_{e} \in \{0, 1\}
            & \forall v \in \matheuler{V},\, 
              \forall e \in \bigcup\rightlap{c \in \matheuler{C}} \matheuler{E}^{v}_{c}
            &
            & \makeref{EVSP2:IntegralityConstraints2}\\
            & ~y^{v}_{e} \in \{0, 1\},\, z^{v}_{e} \in \{0, 1\}
            & \forall v \in \matheuler{V},\, 
              \forall e \in \bigcup\rightlap{s \in \matheuler{S}} \matheuler{E}^{v}_{s}
            &
            & \makeref{EVSP2:IntegralityConstraints4}   
\end{alignat}
\endgroup

The objective function~\eqrefsmallsize{EVSP2:ObjectiveFunction}, the constraint sets~\eqrefsmallsize{EVSP2:DemandConstraints} to~\eqrefsmallsize{EVSP2:InitialDistribution3} and the variable domains~\eqrefsmallsize{EVSP2:NonNegativeConstraints}
to~\eqrefsmallsize{EVSP2:IntegralityConstraints4} are analogous to the first formulation.

\subsubsection*{Formulation \RemoveSpaces{\ref*{formulation:second_splitted}}}

Let us now denote by~\RemoveSpaces{\makerefSb{formulation:second_splitted}} the formulation we obtain by projecting the same set of constraints of~\RemoveSpaces{\ref{formulation:first_splitted}} over the space-time network for each vehicle~\mbox{\(v \in \matheuler{V}\)} created as follows.

For each~\mbox{\(s \in \matheuler{S}\)}, we split station~\(s\) into two separate stations: on the one hand, a station with no charging facilities, with capacity~\mbox{\(\textsf{C}_s - \textsf{R}_s\)}; on the other hand, a station where all parking spaces are equipped with charging facilities, with capacity~\(\textsf{R}_s\).~For each demand in~\(\matheuler{D}_c\), for all~\mbox{\(c \in \matheuler{C}\)}, we split it into four demands, where each represents a combination of pick up at the corresponding station either without or with charging facilities and drop off at the corresponding station either without or with charging facilities.

To further reduce the number of nodes and connecting arcs, we apply a preprocessing step on demands in which the pick-up or drop-off stations have no charging facilities.~Consider the set~\mbox{\(\matheuler{S}^\prime \subseteq \matheuler{S}\)} of all stations without charging facilities, and for each station~\mbox{\(s \in \matheuler{S}^\prime\)}, let~\mbox{\(\outgoing{\matheuler{T}}_{s} = \{t_{i} : \exists (s, t_{i}, \incoming{s}, t_{j}) \in \matheuler{D}\}\)} and~\mbox{\(\incoming{\matheuler{T}}_{s} = \{t_{j} : \exists (\outgoing{s}, t_{i}, s, t_{j}) \in \matheuler{D}\}\)} denote the sets of time instants for the departure and return times, respectively.~This preprocessing step proceeds as follows.

For each demand~\mbox{\((\outgoing{s}, t_{i}, \incoming{s}, t_{j}) \in \matheuler{D}\)}, such that the station~\mbox{\(\outgoing{s} \in \matheuler{S}^\prime\)}, replace the departure time~\(t_{i}\) by~\mbox{\(t_{i^{\prime}} = \min\{t_{i^{\prime}} \in \outgoing{\matheuler{T}}_{s} : t_{i^{\prime}} \geq t_{j^{\prime}}\}\)}, where~\mbox{\(t_{j^{\prime}} = \max\{0, t_{j^{\prime}} \in \incoming{\matheuler{T}}_{s} : t_{j^{\prime}} \leq t_{i}\}\)}.~Note that~\(t_{j^{\prime}}\) is the last return time before the departure time~\(t_{i}\), and that~\(t_{i^{\prime}}\) is the first departure time after the return time~\(t_{j^{\prime}}\).~Replacing~\(t_{i}\) by~\(t_{i^{\prime}}\) implies that the demand departure time is artificially anticipating.

Similarly, for each demand~\mbox{\((\outgoing{s}, t_{i}, \incoming{s}, t_{j}) \in \matheuler{D}\)}, and station~\mbox{\(\incoming{s} \in \matheuler{S}^\prime\)}, replace the return time~\(t_{j}\) by~\mbox{\(t_{i^{\prime}} = \min\{t_{m}, t_{i^{\prime}} \in \outgoing{\matheuler{T}}_{s} : t_{i^{\prime}} \geq t_{j}\}\)}.~Here,~\(t_{i^{\prime}}\) is the first departure time after the return time~\(t_{j}\), and replacing~\(t_{j}\) by~\(t_{i^{\prime}}\) means artificially postponing the demand return time.

\begin{figure}[H]
  \centering
  \resizebox{\textwidth}{!}{%
    \input{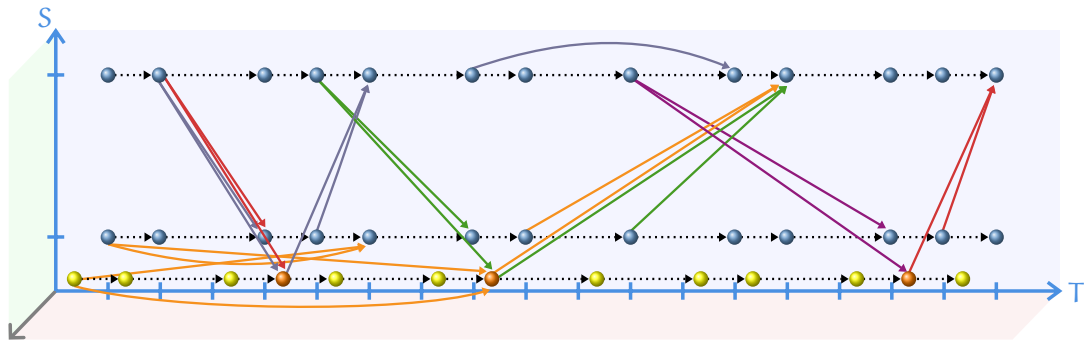}
  }
  \caption{Example of a space-time network modified by splitting the stations containing parking spaces equipped and unequipped with charging facilities and by shifting the appropriate demands.~The blue nodes stand for stations with charging facilities, and the other nodes stand for stations without charging facilities.~In addition, the orange nodes aggregate the shifted demand arcs.}
  \label{fig:example-3D}
\end{figure}

As an illustration, \Cref{fig:example-3D} gives an example for the instance transformation just described, using the network shown in~\Cref{fig:network} as input.~In this illustration, only one of the two stations has all the parking spaces with a charging facility.

Lastly, for all vehicles in~\(\matheuler{V}\), we create their space-time networks using the sets of achievable demands as described before.~Also, all nodes~(and arcs) associated with stations whose capacity has been defined to be zero are removed from the networks.

In particular, in the formulation~\RemoveSpaces{\ref{formulation:second_splitted}}, the constraint set analogous to~\eqrefsmallsize{EVSP1splitted:RechargingConstraints} is defined considering the original return time.

\subsubsection*{Formulation strengths}

For any formulation P, we denote by~$\matheuler{F}({\textnormal{P})}$ the set of feasible solutions of~P.~An immediate consequence is that~$\matheuler{F}$(\RemoveSpaces{\ref*{formulation:second}$_{\textnormal{L}}$}) is a subset of~$\matheuler{F}$(\RemoveSpaces{\ref*{formulation:first}$_{\textnormal{L}}$}) ---~there is an implicit mapping of the variable value from one formulation to the other.~We complement this result by showing that~$\matheuler{F}$(\RemoveSpaces{\ref*{formulation:second}$_{\textnormal{L}}$}) is actually a proper subset of~$\matheuler{F}$(\RemoveSpaces{\ref*{formulation:first}$_{\textnormal{L}}$}).~Consider an instance with two stations, denoted here by \(0\) and \(1\), in which there is initially only one vehicle parked at station~\(0\), and two customers, each with a demand: one demand from station~\(0\) to~\(1\), and immediately after another from station~\(1\) to~\(0\).~Also, consider that the energy required for each demand is strictly less than the energy of the vehicle and that only the first demand is achievable.~For this instance, there is an optimal fractional solution for~\RemoveSpaces{\ref*{formulation:first}$_{\textnormal{L}}$} with a positive flow value through the arc associated with the second demand, but this solution is an infeasible solution for~\RemoveSpaces{\ref*{formulation:second}$_{\textnormal{L}}$}.~This leads to the following result.

\begin{corollary}
  {\normalfont opt(\RemoveSpaces{\ref*{formulation:second}$_{\textnormal{L}}$})} \(\leq\)  {\normalfont opt(\RemoveSpaces{\ref*{formulation:first}$_{\textnormal{L}}$})}.
\end{corollary}

On the other hand, our experiments show that the sets~$\matheuler{F}$(\RemoveSpaces{\ref*{formulation:second}$_{\textnormal{L}}$}) and~$\matheuler{F}$(\RemoveSpaces{\ref*{formulation:first_splitted}$_{\textnormal{L}}$}) are not subsets of each other.

To show that the formulation~\RemoveSpaces{\ref{formulation:second_splitted}} is stronger than the other formulations previously presented for the~\shortname, let us consider the following sequence of propositions.

\begin{proposition}
  Given any feasible solution of \RemoveSpaces{\textnormal{\ref*{formulation:second_splitted}$_{\textnormal{L}}$}}, one can find a feasible solution of~\RemoveSpaces{\textnormal{\ref*{formulation:second}$_{\textnormal{L}}$}} with the same objective value.
\end{proposition}

{\noindent {\bf Proof.}~Consider a feasible solution for~\RemoveSpaces{\ref*{formulation:second_splitted}$_{\textnormal{L}}$}.~To provide an appropriate mapping, consider the following iterative approach.~In the first step, for each variable~\(x^{v}_{e}\), make its value equal to the sum of the values of all variables in the solution, so that they are associated with the same demand for the variable~\(x^{v}_{e}\).~Note that by completing this first step, we ensure that it has the same objective value.

To complete the mapping, make the value of each variable~\(y^{v}_{\outgoing{e}\spacetimezeroskips{s}{0}}\) and~\(z^{v}_{\outgoing{e}\spacetimezeroskips{s}{i}}\), for all~\mbox{\(i \in \matheuler{I}^{v}\)}, equal to the value of the corresponding variable in the solution for~\RemoveSpaces{\ref*{formulation:second_splitted}$_{\textnormal{L}}$}.~Then, iteratively, make the value of each variable~\(y^{v}_{\outgoing{e}\spacetimezeroskips{s}{i}}\) equal to
\begingroup
\small
\begin{equation*}
  \sum\rightlap{e \in \incoming{\matheuler{E}}\spacetimezeroskips{s}{i}} x^{v}_{e} +
  y^{v}_{\incoming{e}\spacetimezeroskips{s}{i}} + z^{v}_{\incoming{e}\spacetimezeroskips{s}{i}} -
  \sum\rightlap{e \in \outgoing{\matheuler{E}}\spacetimezeroskips{s}{i}} x^{v}_{e} -
  z^{v}_{\outgoing{e}\spacetimezeroskips{s}{i}}
\end{equation*}
\endgroup

\noindent which is derived from~\eqrefsmallsize{EVSP2:FlowConservationEquation}, for all~\mbox{\(i \in \matheuler{I}^{v}\!\setminus\!\{0\}\)}.~Note that when iterating in increasing order of the time instants, the variable~\(y^{v}_{\incoming{e}\spacetimezeroskips{s}{i}}\) has its value previously assigned.~The values of the variables~\(\ell^{v}_{i}\) and~\(w_{c}\) are made equal to the values of the corresponding variables in the solution for~\RemoveSpaces{\ref*{formulation:second_splitted}$_{\textnormal{L}}$}.

From here, the rest of the proof follows in the same way as~\Cref{proposition:EVSP1L_to_EVSP1-SL}.~\hfill \(\square\)
}

\begin{proposition}\label{proposition:EVSP2-SL_is_stronger}
  There are instances for which~{\normalfont opt(\RemoveSpaces{\ref*{formulation:second_splitted}$_{\textnormal{L}}$})} is smaller than~{\normalfont opt(\RemoveSpaces{\ref*{formulation:second}$_{\textnormal{L}}$})}.
\end{proposition}

{\noindent {\bf Proof.}~See the~\hyperref[appendix]{Appendix}.}~\hfill \(\square\)

By transitivity, we obtain the following result.

\begin{corollary}
  {\normalfont opt(\RemoveSpaces{\ref*{formulation:second_splitted}$_{\textnormal{L}}$})} \(\leq\)  {\normalfont opt(\RemoveSpaces{\ref*{formulation:first}$_{\textnormal{L}}$})}.
\end{corollary}

Clearly, the set~$\matheuler{F}$(\RemoveSpaces{\ref*{formulation:second_splitted}$_{\textnormal{L}}$}) of feasible solutions is a proper subset of~$\matheuler{F}$(\RemoveSpaces{\ref*{formulation:first_splitted}$_{\textnormal{L}}$}), and we can also obtain the following result.

\begin{corollary}
  {\normalfont opt(\RemoveSpaces{\ref*{formulation:second_splitted}$_{\textnormal{L}}$})} \(\leq\)  {\normalfont opt(\RemoveSpaces{\ref*{formulation:first_splitted}$_{\textnormal{L}}$})}.
\end{corollary}

\begin{figure}[H]
  \centering
  \resizebox{.4\textwidth}{!}{%
    \tikzset{every picture/.style={line width=0.75pt}} 

\begin{tikzpicture}[x=0.75pt,y=0.75pt,yscale=-1,xscale=1]

\draw (160.6,50.73) -- (249.4,50.73) -- (249.4,89.4) -- (160.6,89.4) -- cycle;
\draw (81.2,120.33) -- (170,120.33) -- (170,159) -- (81.2,159) -- cycle;
\draw (160.8,190.93) -- (249.6,190.93) -- (249.6,229.6) -- (160.8,229.6) -- cycle;
\draw (240.4,121.07) -- (329.2,121.07) -- (329.2,159.73) -- (240.4,159.73) -- cycle;

\draw [arrows={-angle 45}] (206.2-7,190.8) -- (127.28-6,161.89);
\draw [arrows={-angle 45}] (206.2+3,190.8) -- (284.32+6,162.29);
\draw [arrows={-angle 45}] (122.6,120.4) -- (200.72,91.89);
\draw [arrows={-angle 45}] (285.8,121) -- (206.88+2,92.09);
\draw [dash pattern={on 4.5pt off 4.5pt}] (172,140) -- (242.6,140);

\draw (284.71,140.4) node   [align=left] {\begin{minipage}[lt]{60.26pt}\setlength\topsep{0pt}
   \centering\RemoveSpaces{\ref*{formulation:second}$_{\textnormal{L}}$}
\end{minipage}};
\draw (205.11,210.27) node   [align=left] {\begin{minipage}[lt]{60.26pt}\setlength\topsep{0pt}
   \centering\RemoveSpaces{\ref*{formulation:first}$_{\textnormal{L}}$}
\end{minipage}};
\draw (125.51,139.67) node   [align=left] {\begin{minipage}[lt]{60.26pt}\setlength\topsep{0pt}
   \centering\RemoveSpaces{\ref*{formulation:first_splitted}$_{\textnormal{L}}$}
\end{minipage}};
\draw (204.91,70.07) node   [align=left] {\begin{minipage}[lt]{60.26pt}\setlength\topsep{0pt}
   \centering\RemoveSpaces{\ref*{formulation:second_splitted}$_{\textnormal{L}}$}
\end{minipage}};

\end{tikzpicture}
  }
  \caption{Dominance relationships between {\shortname} formulations.~The arrows lead from dominated formulations to dominating ones.~Formulations without a dominance relationship are shown at the same level, with a dashed line.}
  \label{fig:dominance_relationships}
\end{figure}
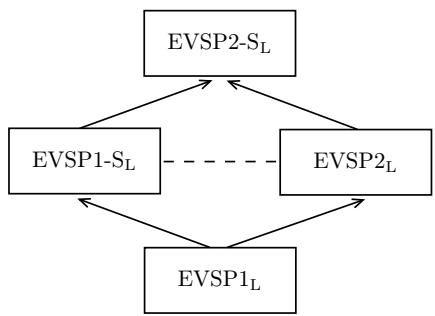

\Cref{fig:dominance_relationships} illustrates the dominance relationships of the linear programming relaxations of the formulations. 

\section{Complexity}
\label{sec:complexity}

In this section, we discuss the computational complexity of our problem.~The~\mbox{NP-hardness} of the~\shortname can be proved using the same arguments used to prove the hardness of the TFC problem~\citep{Bohmova/2016}.~The proof is based on a special case where the driving time is constant, and there is only one vehicle.~In this special case, maximizing the number of customers served is equivalent to maximizing the total vehicle rental time.~However, the authors show that the~TFC problem can be solved in polynomial time when each customer has only one demand.~An interesting question is whether the~\shortname is also solvable in polynomial time when each customer has only one demand. 

We prove that the~\shortname is~\mbox{NP-hard} even if each customer has only one demand, by a reduction from the one-dimensional bin packing problem~(BPP) described as follows.

Let~\mbox{\(\matheuler{N} = \{1, \dots, n\}\)} be a set of one-dimensional items, where each item~\mbox{\(i \in \matheuler{N}\)} has a size~\mbox{\(\ell_{i} \in (0, 1]\)}. Given bins of capacity~\(1\), the objective of the BPP is to pack the~\(n\) items into the smallest number of bins, such that no bin capacity is exceeded.~Equivalently, this problem aims to find a partition~\mbox{\(\{ \matheuler{N}_{1}, \dots, \matheuler{N}_{k} \}\)} of~\(\matheuler{N}\), such that the sum of the sizes of the items in each part is at most~\(1\) and that~\(k\) is minimized.

\begin{figure}[H]
  \centering
  \resizebox{\textwidth}{!}{%
    \tikzset{every picture/.style={line width=0.75pt}}

\begin{tikzpicture}[x=0.75pt, y=0.75pt, yscale=-1, xscale=0.96]

   \draw  [draw opacity=0][fill={rgb, 255:red, 185; green, 215; blue, 252}  ,fill opacity=1] (24.67,112.33) -- (561.17,112.33) -- (561.17,102.33) -- (584,122.33) -- (561.17,142.33) -- (561.17,132.33) -- (24.67,132.33) -- cycle;

   \draw [->,>={Triangle[scale=.6]},line width=.4mm,dash pattern={on 5.63pt off 4.5pt}] (51,112.67) .. controls (76.35,80.33) and (97.9,77.94) .. (127.07,109.97);
   
   \draw [->,>={Triangle[scale=.6]},line width=.35mm] (130.67,131.67) .. controls (162.35,164.82) and (182.2,163.57) .. (207.45,133.38);
   \draw [->,>={Triangle[scale=.6]},line width=.35mm] (291.33,131) .. controls (323.02,164.16) and (342.86,162.9) .. (368.11,132.72);
   \draw [->,>={Triangle[scale=.6]},line width=.35mm] (451.33,131) .. controls (483.02,164.16) and (502.86,162.9) .. (528.11,132.72);
   
   \draw (41.84,116.44) node [anchor=north west][inner sep=0.75pt][align=left] {$\displaystyle t_{1}$};
   \draw (121.84,116.44) node [anchor=north west][inner sep=0.75pt][align=left] {$\displaystyle t_{2}$};
   \draw (201.84,116.44) node [anchor=north west][inner sep=0.75pt][align=left] {$\displaystyle t_{3}$};
   \draw (282.25,116.44) node [anchor=north west][inner sep=0.75pt][align=left] {$\displaystyle t_{4}$};
   \draw (362.25,116.44) node [anchor=north west][inner sep=0.75pt][align=left] {$\displaystyle t_{5}$};
   \draw (442.25,116.44) node [anchor=north west][inner sep=0.75pt][align=left] {$\displaystyle t_{2n}$};
   \draw (521.92,116.44) node [anchor=north west][inner sep=0.75pt][align=left] {$\displaystyle t_{2n+1}$};
   \draw (397,121.67) node [anchor=north west][inner sep=0.75pt][align=left] {$\displaystyle \dotsc $};

\end{tikzpicture}
  }
  \caption{Placement of demands for the~BPP instance transformed into an instance of the~\shortname.~The demands of item customers are displayed as full arrows and for the demands of dummy customers there is only one dashed arrow for simplicity.}
  \label{fig:BPPtoEVSP}
\end{figure}
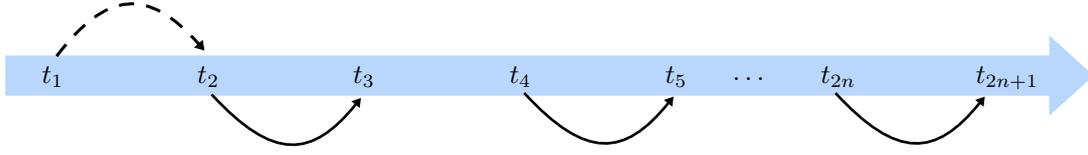

Consider an arbitrary instance of the~BPP, denoted by~\(\matheuler{I}_{\textsf{BPP}}\).~From~\(\matheuler{I}_{\textsf{BPP}}\), we construct in polynomial time an instance of the~\shortname, denoted by~\(\matheuler{I}_{\textsf{\shortname}}\), in which there are~\(2n\) customers and each customer has a single demand.~In~\(\matheuler{I}_{\textsf{\shortname}}\), there is only one station, denoted by~\(s\), with capacity~\(n\) and no charging facilities.~There are~\(n\) fully charged vehicles initially located at~\(s\) with battery capacity equal to~\(1\)~(\ie, equal to the capacity of the bin).~Each item~\mbox{\(i \in \matheuler{N}\)} is associated with an~\emph{item customer} with a demand~\mbox{\((s, 2i, s, 2i\!+\!1, \ell_{i})\)}.~There are also~\(n\)~\emph{dummy customers}, each with a demand~\mbox{\((s, 1, s, 2, \textsf{L})\)}.~See~\Cref{fig:BPPtoEVSP} for an exact configuration.

Note that each customer contributes with a rental time of one unit in the objective function.~Thus, maximizing the total vehicle rental time is equivalent to maximizing the number of customers served.

We now describe the properties of the instance~\(\matheuler{I}_{\textsf{\shortname}}\) constructed from~\(\matheuler{I}_{\textsf{BPP}}\), defined as propositions.~Let~\(f\) be a polynomial time reduction from~BPP to~\shortname as described above, such that~\mbox{\(\matheuler{I}_{\textsf{\shortname}} = f(\matheuler{I}_{\textsf{BPP}})\)}.

\begin{proposition}
  \label{proposition:k_vehicles}
  Let~\(\matheuler{I}_{\normalfont \textsf{BPP}}\) be an instance of the~{\normalfont \textsf{BPP}}, and consider a solution for~\(\matheuler{I}_{\normalfont \textsf{BPP}}\) where~\(k\) bins are used.~If~\mbox{\(\matheuler{I}_{\normalfont \textsf{\shortname}} = f(\matheuler{I}_{\normalfont \textsf{BPP}})\)}, then there is a solution for~\(\matheuler{I}_{\normalfont \textsf{\shortname}}\)\! such that~\(k\) vehicles are used to serve all item customers.
\end{proposition}

{\noindent {\bf Proof.}~Consider a solution for~\(\matheuler{I}_{\textsf{BPP}}\) represented by a partition~\mbox{\(\{\matheuler{N}_{1}, \dots, \matheuler{N}_{k}\}\)} of~\(\matheuler{N}\).~By the construction of~\(\matheuler{I}_{\textsf{\shortname}}\), all item customers for which there is an item associated in~\(\matheuler{N}_{i}\), for any~\mbox{\(1 \leq i \leq k\)}, can be served by the same vehicle, consuming at most the battery capacity.~\hfill \(\square\)
}


\begin{proposition}
  \label{proposition:opt_value}
  Let~{\normalfont opt}\((\matheuler{I}_{\normalfont \textsf{\shortname}})\)\! and~{\normalfont opt}\((\matheuler{I}_{\normalfont \textsf{BPP}})\)\! denote, respectively, the optimal value for the instances~\(\matheuler{I}_{\normalfont \textsf{\shortname}}\)\! and~\(\matheuler{I}_{\normalfont \textsf{BPP}}\), where~\mbox{\(\matheuler{I}_{\normalfont \textsf{\shortname}} = f(\matheuler{I}_{\normalfont \textsf{BPP}})\)}.~Then~\mbox{{\normalfont opt}\((\matheuler{I}_{\normalfont \textsf{\shortname}}) = 2n\, -\)\,{\normalfont opt}\((\matheuler{I}_{\normalfont \textsf{BPP}})\)}.
\end{proposition}

{\noindent {\bf Proof.}~Let~\mbox{\(k =\) opt\((\matheuler{I}_{\textsf{BPP}})\)}.~By~\Cref{proposition:k_vehicles}, we know that there exists a solution for~\(\matheuler{I}_{\textsf{\shortname}}\) in which all item customers are served by exactly~\(k\) vehicles.~Thus, there are~\mbox{\(n - k\)} vehicles that are not used in this solution.~Since each of these vehicles can only serve one dummy customer, we can then state that an optimal solution for~\(\matheuler{I}_{\textsf{\shortname}}\) must serve at least~\mbox{\(2n - k\)} customers.~That is,~\mbox{opt\((\matheuler{I}_{\textsf{\shortname}}) \geq 2n\, -\)\,opt\((\matheuler{I}_{\textsf{BPP}})\)}.

Now we show that~\mbox{opt\((\matheuler{I}_{\textsf{\shortname}}) \leq 2n\, -\)\,opt\((\matheuler{I}_{\textsf{BPP}})\)}.~Consider an optimal solution for~\(\matheuler{I}_{\textsf{\shortname}}\) where the number of dummy customers served is as minimum as possible.~Note that in any optimal solution for~\(\matheuler{I}_{\textsf{\shortname}}\), all~\(n\) vehicles are used, since an idle vehicle could always be used to serve at least one customer.~Suppose that there is at least one dummy customer in this optimal solution.~Thus, all item customers are served, since if there is any item customer unserved, then the number of dummy customers served is not minimum.~Therefore, if there are~\(k\) vehicles to serve all~\(n\) item customers, and there are~\mbox{\(n - k\)} vehicles to serve~\mbox{\(n - k\)} dummy customers, then~\mbox{opt\((\matheuler{I}_{\textsf{\shortname}}) = 2n - k\)}.

Note that all item customers served by the same vehicle stand for a bin in a solution for~\(\matheuler{I}_{\textsf{BPP}}\), therefore,~\mbox{opt\((\matheuler{I}_{\textsf{BPP}}) \leq k = 2n -\)\,opt\((\matheuler{I}_{\textsf{\shortname}})\)}.~\hfill \(\square\)
} 

From an optimal solution for~\mbox{\(\matheuler{I}_{\textsf{\shortname}} = f(\matheuler{I}_{\textsf{BPP}})\)}, we can construct an optimal solution for~\(\matheuler{I}_{\textsf{BPP}}\) in polynomial time, as follows: for each vehicle used to serve item customers, we pack all items associated with customers served by the vehicle into a single bin; and we pack each item associated with an unserved item customer, if any, into its own bin.~Then, the next proposition holds that~\(\matheuler{I}_{\textsf{BPP}}\) is an optimal solution.

\begin{proposition}
  Let~\(\matheuler{S}_{\normalfont \textsf{\shortname}}\)\! denote an optimal solution for~\mbox{\(\matheuler{I}_{\normalfont \textsf{\shortname}} = f(\matheuler{I}_{\textsf{BPP}})\)}.~If a solution for~\(\matheuler{I}_{\normalfont \textsf{BPP}}\), denoted by~\(\matheuler{S}_{\normalfont \textsf{BPP}}\), is constructed from~\(\matheuler{S}_{\normalfont \textsf{\shortname}}\) as described above, then~\(\matheuler{S}_{\normalfont \textsf{BPP}}\)\! is an optimal solution.
\end{proposition}

{\noindent {\bf Proof.}~Suppose, without loss of generality, that in~\(\matheuler{S}_{\textsf{\shortname}}\) all item customers are served.~From this, the proof follows directly from~\Cref{proposition:opt_value}.~\hfill \(\square\)
}

Finally, as a corollary, we obtain the following hardness result for our problem using the previous propositions.

\begin{corollary}
  \label{proposition:csp}
  The~{\normalfont \shortname}\! is~\mbox{{\normalfont NP}-hard} even if every customer has only one demand and when there is only one station.
\end{corollary}



\section{Computational experiments} 
\label{sec:experiments}

We dedicate this section to evaluate the efficiency of the proposed mixed-integer linear programming formulations for solving the~\shortname using two benchmark sets.~The first set was obtained by adapting the random instance generation procedure used by~\cite{Xu/2019}, while the second set was obtained using data from the~VAMO Fortaleza system~\citep{Fortaleza/2018}.~Both sets of instances are available at~\href{https://gitlab.com/welverton/evsp}{https:{\small/\!/}gitlab.com{\small/}welverton{\small/}evsp}.

In what follows, we drop the integrality constraints only for the variables for which their integrality is ensured by the other constraints in all optimal solutions as described earlier.~Although, the commercial solver might classify some of those relaxed variables as implicit integer variables\footnote{A continuous variable is an~\emph{implicit integer variable} if it can only take integral values when all integer variables take integral values.~Implicit integer variables do not have to be considered for branching decisions.} during the presolve phase
.

We provide an initial solution as a warm-start for each execution.~The initial solution is obtained by a fast construction heuristic, that first sorts the customers in descending order of their total rental time divided by their number of demands.~Iteratively, each customer is assigned to a vehicle in that order if and only if the vehicle can fulfill all customer demands without violating any constraints.~Vehicle recharging is not taken into account when creating the initial solution.

We assume that all vehicles have a fully charged battery at the beginning of the daily planning and that the vehicles are initially distributed at each station such that the number of empty parking spaces equipped with charging facilities is maximum.

All computational experiments were performed on an~\mbox{Ubuntu 18.04 64-bit} operating system using an~Intel(R) Xeon(R) E5-2630 v4 running at~\mbox{2.2 GHz}, \(10\)~cores, \mbox{and \(64\) GB RAM}.~Our implementation was done in~{\CC} and uses~Gurobi Optimizer \mbox{version 9.0.3} as the mixed-integer linear programming solver, compiled with~GCC version 7 using the optimization flag~-O3.~We used Gurobi Optimizer with the default settings and restricted the computation to a single thread.~A time limit of one hour was set for each execution.

\subsection{General performance analysis}

To create the random instances, we assume an average speed \mbox{of \(30\) km} per hour, denoted by~\(\vartheta\).~Thus the minimum demand rental time is~\mbox{\(dis(\outgoing{s}, \incoming{s})/ \vartheta\)} hours, where~\mbox{\(dis(\outgoing{s}, \incoming{s})\)} denotes the Euclidean distance between the pick-up and the drop-off stations of the demand.~In this case, the pick-up and the drop-off stations of each demand are different stations.~The battery capacity of each vehicle is~\mbox{\(40\) kilowatt-hours (kWh)} and the discharge rate is~\mbox{\(5\) kWh}~(denoted as~\(\lambda\)), which is roughly equivalent to a~Nissan Leaf \(40\)~kWh~\citep{Nissan/2019}.~Assuming~\(4\) hours as the maximum time to fully charge the battery, the energy supplied by the charging facility \mbox{is \(10\) kWh} per hour.~We assume an operating period from~\mbox{5:00} to~\mbox{23:00} per day, with a morning peak period from~\mbox{7:00} to~\mbox{10:00} and an evening peak period from~\mbox{16:00} to~\mbox{20:00}, so that the driving demands during off-peak hours are low.

Our first set of instances can be divided into groups of~\mbox{\(30, 40, \dots, 120\)} customers, with each group containing instances ranging from~\(3\) to~\(5\) stations.~All instances are created as described below.

\noindent \hspace*{-.5mm}\emph{\textbf{Grid instances}.}~At the first step,~\(\cardinality{\matheuler{S}}\) stations are chosen uniformly at random from a~\(50\)~km by~\(50\)~km grid.~For each station~\mbox{\(s \in \matheuler{S}\)}, the capacity~\(\textsf{C}_{s}\) is given by a random integer between~\(10\) and~\(20\), while the number~\(\textsf{R}_{s}\) of charging facilities is  an integer chosen randomly between~\(20\%\) and~\(50\%\) of the station capacity (\ie, limited availability of charging infrastructure).~Similarly, the number of vehicles initially located at a station is a random integer between~\(30\%\) and~\(75\%\) of its capacity.~Furthermore,~\(\cardinality{\matheuler{C}}\) customers have their demands defined.~Each set of demands contains at most four demands, in which each possible size has a predefined probability of being produced.~The weights used are~\(4\), \(6\), \(2\), and~\(1\), for sizes~\(1\), \(2\), \(3\), and~\(4\), respectively (each weight unit provides a~\(1/13\) probability).~For each demand~\mbox{\((\outgoing{s}, t_{i}, \incoming{s}, t_{j}, \varepsilon) \in \matheuler{D}_{c}\)}, the pick-up and drop-off stations are randomly chosen from the generated stations.~The departure time is defined by first randomly selecting one of~\(18\) one-hour intervals of the operating period, in which each interval has a predefined probability of being selected.~We used weights~\mbox{\((4, 7, 3)\)} for the morning peak intervals, weights~\mbox{\((2, 5, 7, 3)\)} for the evening peak intervals, and weight~\(1\) for each off-peak interval.~The departure time is then randomly chosen from the set of time instants~(with a granularity of ﬁve minutes) associated with the selected one-hour interval.~The arrival time is thus chosen as a random integer from the set~\mbox{\(\{dis(\outgoing{s}, \incoming{s})/\vartheta + 5 \text{ min}, dis(\outgoing{s}, \incoming{s})/\vartheta + 10 \text{ min}, \dots,\) \(dis(\outgoing{s}, \incoming{s})/\vartheta + 30 \text{ min}\}\)}.~The energy required for the rental period is randomly chosen from the interval~\mbox{\([dis(\outgoing{s}, \incoming{s})/\vartheta,\, t_{j}-t_{i})\times\lambda\)}.

\begin{figure}[H]
  \centering
  \begin{tikzpicture}[scale=0.83]
    \begin{axis}[
          xlabel={Time (minutes)},
          ylabel={Number of solved instances (\%)},
          ylabel style={at={(0.03,0.5)}},
          legend entries={
            \RemoveSpaces{\ref*{formulation:first}},
            \RemoveSpaces{\ref*{formulation:first_splitted}},
            \RemoveSpaces{\ref*{formulation:second}},
            \RemoveSpaces{\ref*{formulation:second_splitted}},
          },
          legend style={at={(.97,.92)},anchor=north east},
          legend cell align=left,
          x=0.1085cm, y=0.06cm,
          scale only axis, 
          ymin=0,ymax=104,xmin=0,xmax=60,
          enlargelimits=false,
      ]
      \addplot[mark=none, Red, line width=1pt] table{data/Grid/performance/EVSP1.dat};
      \addplot[mark=none, MidnightBlue, line width=1pt] table{data/Grid/performance/EVSP1-S.dat};
      \addplot[mark=none, Green, line width=1pt] table{data/Grid/performance/EVSP2.dat};
      \addplot[mark=none, Goldenrod, line width=1pt] table{data/Grid/performance/EVSP2-S.dat};
    \end{axis}
  \end{tikzpicture}
  \qquad
  \begin{tikzpicture}[scale=0.83]
    \begin{axis}[
          xlabel={Optimality gap (\%)},
          ylabel={Number of instances (\%)},
          ylabel style={at={(0.03,0.5)}},
          legend entries={
            \RemoveSpaces{\ref*{formulation:first}},
            \RemoveSpaces{\ref*{formulation:first_splitted}},
            \RemoveSpaces{\ref*{formulation:second}},
            \RemoveSpaces{\ref*{formulation:second_splitted}},
          },
          legend style={at={(.97,.92)},anchor=north east},
          legend cell align=left,
          x=0.0186cm, y=0.06cm,
          scale only axis, 
          ymin=0,ymax=104,xmin=0,xmax=350,
          enlargelimits=false,
        ]
        \addplot[mark=none, Red, line width=1pt] table{data/Grid/gap/EVSP1.dat};
        \addplot[mark=none, MidnightBlue, line width=1pt] table{data/Grid/gap/EVSP1-S.dat};
        \addplot[mark=none, Green, line width=1pt] table{data/Grid/gap/EVSP2.dat};
        \addplot[mark=none, Goldenrod, line width=1pt] table{data/Grid/gap/EVSP2-S.dat};
    \end{axis}
  \end{tikzpicture}
  \caption{Percentage of optimally solved instances and relative optimality gaps within the one-hour time limit.~Gap values were calculated as~\mbox{(\hspace*{-.1mm}(\hspace*{-.1mm}\textsf{UB} - \textsf{LB})/\textsf{LB})\(\times 100\))}, where \textsf{UB} and~\textsf{LB} are the upper and lower bound values achieved by the solver, respectively.}
  \label{fig:grid:performace_profile}
\end{figure}
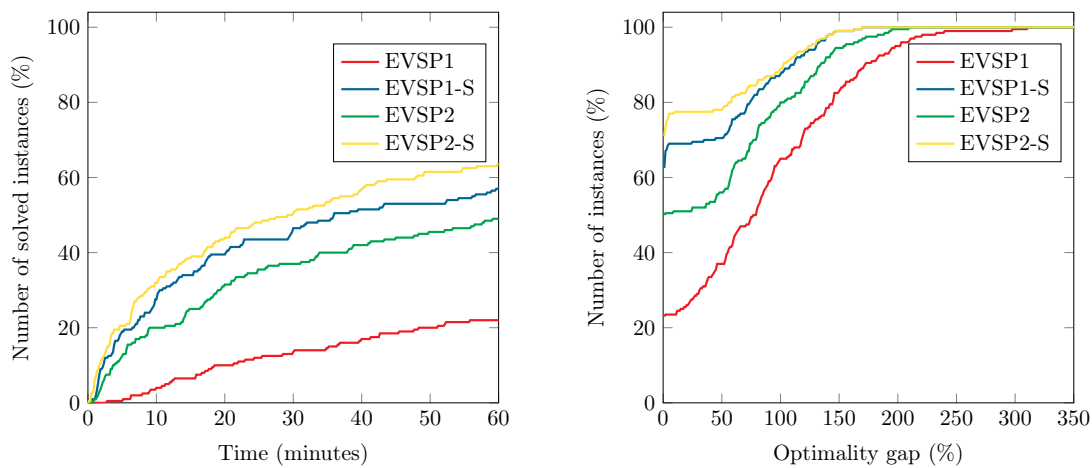

Now we present the computational results for these instances.~\Cref{fig:grid:performace_profile} shows the percentage of solved instances and the remaining optimality gaps.~From these plots, we can see that~\RemoveSpaces{\ref{formulation:first_splitted}} and~\RemoveSpaces{\ref{formulation:second_splitted}} were significantly faster than the other two formulations and also obtained a better gap.~Furthermore, these models were the only ones that could solve all instances with up to~\(50\) customers optimally.

\begin{figure}[H]
  \centering
  \begin{center}
  \begin{tikzpicture}[scale=0.83]
    \begin{axis}[box plot width=1mm, xtick distance=1, 
                  ylabel={Number of binary variables},
                  ylabel style={at={(0.03,0.5)}},
                  xticklabels={,,
                               {\small \ref*{formulation:first}},
                               {\small \ref*{formulation:first_splitted}},
                               {\small \ref*{formulation:second}},
                               {\small \ref*{formulation:second_splitted}}
                              }
                  ]
      \boxplot [
        forget plot, gray, xshift=-.25cm,
        box plot whisker bottom index=1,
        box plot whisker top index=5,
        box plot box bottom index=2,
        box plot box top index=4,
        box plot median index=3
      ] {binary:original.dat}
      \boxplot [
        forget plot, red, xshift=.0cm,
        box plot whisker bottom index=1,
        box plot whisker top index=5,
        box plot box bottom index=2,
        box plot box top index=4,
        box plot median index=3
      ] {binary:relax.dat}
      \boxplot [
        forget plot, MidnightBlue, xshift=.25cm,
        box plot whisker bottom index=1,
        box plot whisker top index=5,
        box plot box bottom index=2,
        box plot box top index=4,
        box plot median index=3
      ] {binary:relax_and_presolved.dat}
    \end{axis}
  \end{tikzpicture}
  \qquad
  \begin{tikzpicture}[scale=0.83]
    \begin{axis}[box plot width=1mm, xtick distance=1, 
                  ylabel={Number of constraints}, 
                  ylabel style={at={(0.03,0.5)}},
                  xticklabels={,,
                               {\small \ref*{formulation:first}},
                               {\small \ref*{formulation:first_splitted}},
                               {\small \ref*{formulation:second}},
                               {\small \ref*{formulation:second_splitted}}
                              }
                ]
      \boxplot [
        forget plot, red, xshift=.0cm,
        box plot whisker bottom index=1,
        box plot whisker top index=5,
        box plot box bottom index=2,
        box plot box top index=4,
        box plot median index=3
      ] {constraints:relax.dat}
      \boxplot [
        forget plot, MidnightBlue, xshift=.25cm,
        box plot whisker bottom index=1,
        box plot whisker top index=5,
        box plot box bottom index=2,
        box plot box top index=4,
        box plot median index=3
      ] {constraints:relax_and_presolved.dat}
    \end{axis}
  \end{tikzpicture}
  \end{center}
  \caption{Distribution of the number of binary variables and constraints on the last group of instances~(\ie, instances with~\(120\) customers).~The red and blue boxplots correspond respectively to the formulations before and after the presolve phase with relaxation of the integrality of the variables that can be relaxed.~In addition, the gray boxplots correspond to the formulations without any relaxation for integrality.}
  \label{fig:grid:dimension}
\end{figure}
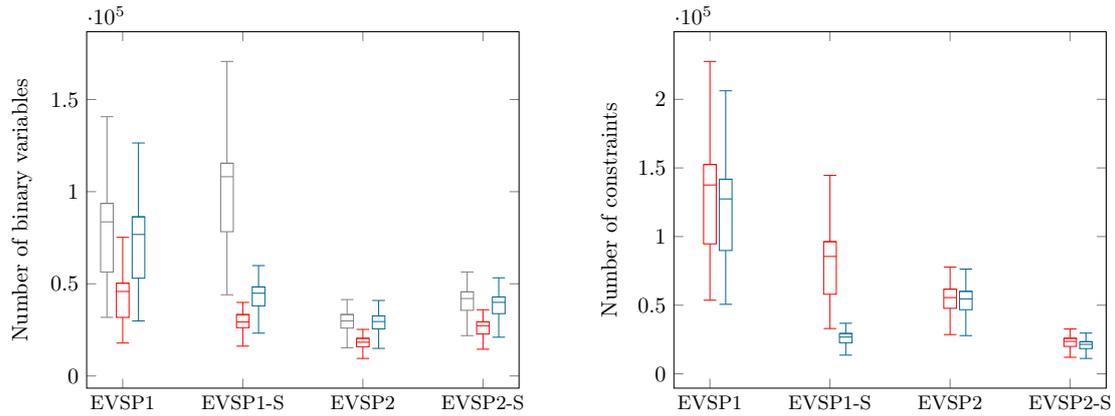

Remember that~\RemoveSpaces{\ref{formulation:first_splitted}} has a relatively large number of binary variables associated with customers' demands compared to the other formulations.~However, the presolve phase can be very effective in reducing the size of the problem in advance and tightening its integer linear program.~\Cref{fig:grid:dimension} shows the boxplot for the number of binary variables~(left plot) and the boxplot for the number of constraints~(right plot).~Each boxplot shows the median~(horizontal line inside the box), while the extremes of the boxes represent the first and third quartiles, respectively, and the whiskers represent the smallest and largest values.~Note that presolve phase of the Gurobi changes some continuous variables (\ie, some relaxed variables to be continuous) to implicit integer variables.

From these boxplots, we can see that~\RemoveSpaces{\ref{formulation:first_splitted}} and~\RemoveSpaces{\ref{formulation:second_splitted}} have approximately the same number of binary variables and constraints after presolve.~\Cref{tab:grid:dimension} provides a comparison between the original size of the integer linear programs and their size after the presolve phase.~Note that the integer linear programs of~\RemoveSpaces{\ref{formulation:first_splitted}} after the presolve phase are larger than~\RemoveSpaces{\ref{formulation:second_splitted}} before the presolve begins~(\ie, the original size of its integer linear programs).~The best results are highlighted with a gray background.

\begin{table}[H]
  \caption{Size comparison of the integer linear programs where the percentages related in each column are relative to the size of~\RemoveSpaces{\ref*{formulation:first_splitted}} before the presolve phase.}
  \label{tab:grid:dimension}
  \setlength{\tabcolsep}{.1em}
  \renewcommand{\arraystretch}{1.1}
  \resizebox{\textwidth}{!}{%
  \begin{tabular}{cccccccccccc}
    \toprule
     & \multicolumn{3}{c}{\normalsize\RemoveSpaces{\ref*{formulation:first_splitted}} after presolve} &  & \multicolumn{3}{c}{\normalsize\RemoveSpaces{\ref*{formulation:second_splitted}} before presolve} &  & \multicolumn{3}{c}{\normalsize\RemoveSpaces{\ref*{formulation:second_splitted}} after presolve} \\ \cline{2-4} \cline{6-8} \cline{10-12} 
    Instance~~ & ~~~Rows~~~    & ~Columns\, & \,Nonzeros\, &~~~~& ~~~Rows~~~    & ~Columns\, & \,Nonzeros\, &~~~~& ~~~Rows~~~    & ~Columns\, & \,Nonzeros\, \\ \midrule \\[-16pt]
    181      & 27.00\% & 40.72\% & 59.58\%  &  & 23.49\% & 38.20\% & 48.94\%  &  & \cellcolor[HTML]{EFEFEF}21.29\% & \cellcolor[HTML]{EFEFEF}36.25\% & \cellcolor[HTML]{EFEFEF}48.69\%  \\
    182      & 47.34\% & 64.07\% & 77.60\%  &  & 41.95\% & 61.06\% & 70.85\%  &  & \cellcolor[HTML]{EFEFEF}39.12\% & \cellcolor[HTML]{EFEFEF}59.07\% & \cellcolor[HTML]{EFEFEF}70.47\%  \\
    183      & 26.99\% & 40.06\% & 58.63\%  &  & 23.46\% & 37.55\% & 47.98\%  &  & \cellcolor[HTML]{EFEFEF}21.46\% & \cellcolor[HTML]{EFEFEF}35.76\% & \cellcolor[HTML]{EFEFEF}47.71\%  \\
    184      & 23.92\% & 35.93\% & 53.75\%  &  & 20.48\% & 33.48\% & 43.16\%  &  & \cellcolor[HTML]{EFEFEF}18.50\% & \cellcolor[HTML]{EFEFEF}31.66\% & \cellcolor[HTML]{EFEFEF}42.80\%  \\
    185      & 33.13\% & 48.22\% & 66.11\%  &  & 29.21\% & 45.56\% & 56.57\%  &  & \cellcolor[HTML]{EFEFEF}26.31\% & \cellcolor[HTML]{EFEFEF}43.18\% & \cellcolor[HTML]{EFEFEF}56.24\%  \\
    186      & 41.26\% & 56.48\% & 71.72\%  &  & 36.45\% & 53.60\% & 63.64\%  &  & \cellcolor[HTML]{EFEFEF}33.92\% & \cellcolor[HTML]{EFEFEF}51.66\% & \cellcolor[HTML]{EFEFEF}63.25\%  \\
    187      & 31.99\% & 46.76\% & 64.56\%  &  & 28.21\% & 44.23\% & 55.01\%  &  & \cellcolor[HTML]{EFEFEF}26.00\% & \cellcolor[HTML]{EFEFEF}42.38\% & \cellcolor[HTML]{EFEFEF}54.72\%  \\
    188      & 42.29\% & 57.42\% & 73.40\%  &  & 37.17\% & 54.15\% & 64.72\%  &  & \cellcolor[HTML]{EFEFEF}34.75\% & \cellcolor[HTML]{EFEFEF}52.26\% & \cellcolor[HTML]{EFEFEF}64.44\%  \\
    189      & 25.30\% & 38.44\% & 57.31\%  &  & 22.40\% & 36.41\% & 46.75\%  &  & \cellcolor[HTML]{EFEFEF}20.04\% & \cellcolor[HTML]{EFEFEF}34.32\% & \cellcolor[HTML]{EFEFEF}46.43\%  \\
    190      & 40.97\% & 56.15\% & 72.22\%  &  & 36.37\% & 53.25\% & 63.74\%  &  & \cellcolor[HTML]{EFEFEF}33.08\% & \cellcolor[HTML]{EFEFEF}50.71\% & \cellcolor[HTML]{EFEFEF}63.37\%  \\
    191      & 32.05\% & 45.78\% & 62.69\%  &  & 27.98\% & 43.04\% & 53.27\%  &  & \cellcolor[HTML]{EFEFEF}25.51\% & \cellcolor[HTML]{EFEFEF}40.98\% & \cellcolor[HTML]{EFEFEF}52.91\%  \\
    192      & 32.02\% & 46.78\% & 64.54\%  &  & 28.06\% & 44.11\% & 54.82\%  &  & \cellcolor[HTML]{EFEFEF}25.45\% & \cellcolor[HTML]{EFEFEF}41.95\% & \cellcolor[HTML]{EFEFEF}54.54\%  \\
    193      & 30.67\% & 45.13\% & 62.79\%  &  & 27.02\% & 42.65\% & 53.21\%  &  & \cellcolor[HTML]{EFEFEF}24.50\% & \cellcolor[HTML]{EFEFEF}40.55\% & \cellcolor[HTML]{EFEFEF}52.90\%  \\
    194      & 30.80\% & 44.38\% & 62.01\%  &  & 27.12\% & 42.15\% & 52.51\%  &  & \cellcolor[HTML]{EFEFEF}24.67\% & \cellcolor[HTML]{EFEFEF}39.78\% & \cellcolor[HTML]{EFEFEF}51.67\%  \\
    195      & 43.95\% & 59.98\% & 75.66\%  &  & 38.50\% & 56.51\% & 67.36\%  &  & \cellcolor[HTML]{EFEFEF}35.42\% & \cellcolor[HTML]{EFEFEF}54.17\% & \cellcolor[HTML]{EFEFEF}67.02\%  \\
    196      & 26.23\% & 39.93\% & 59.02\%  &  & 23.27\% & 37.85\% & 48.50\%  &  & \cellcolor[HTML]{EFEFEF}21.15\% & \cellcolor[HTML]{EFEFEF}35.94\% & \cellcolor[HTML]{EFEFEF}48.16\%  \\
    197      & 43.37\% & 59.47\% & 74.66\%  &  & 38.20\% & 56.26\% & 66.73\%  &  & \cellcolor[HTML]{EFEFEF}35.33\% & \cellcolor[HTML]{EFEFEF}54.08\% & \cellcolor[HTML]{EFEFEF}66.39\%  \\
    198      & 44.09\% & 60.35\% & 75.32\%  &  & 38.91\% & 57.25\% & 67.42\%  &  & \cellcolor[HTML]{EFEFEF}35.59\% & \cellcolor[HTML]{EFEFEF}54.82\% & \cellcolor[HTML]{EFEFEF}67.07\%  \\
    199      & 33.41\% & 48.84\% & 66.09\%  &  & 29.58\% & 46.44\% & 57.12\%  &  & \cellcolor[HTML]{EFEFEF}26.95\% & \cellcolor[HTML]{EFEFEF}44.26\% & \cellcolor[HTML]{EFEFEF}56.72\%  \\
    200      & 34.22\% & 49.77\% & 67.31\%  &  & 30.15\% & 47.08\% & 58.05\%  &  & \cellcolor[HTML]{EFEFEF}27.57\% & \cellcolor[HTML]{EFEFEF}44.95\% & \cellcolor[HTML]{EFEFEF}57.73\%  \\ \toprule
  \multicolumn{12}{l}{\small \emph{Note.} Rows, number of rows; Columns, number of columns; Nonzeros, number of nonzero elements.}
  \end{tabular}}
\end{table}

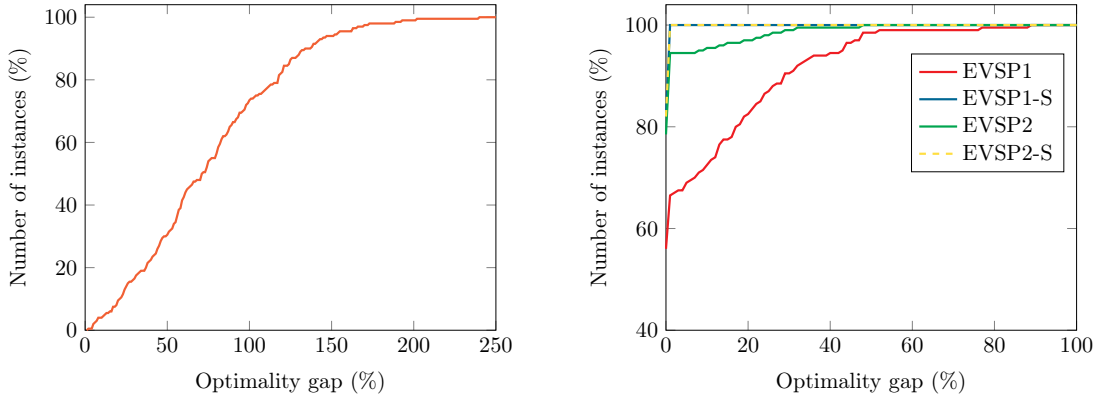
\begin{figure}[H]
  \centering
  \begin{tikzpicture}[scale=0.83]
    \begin{axis}[
          xlabel={Optimality gap (\%)},
          ylabel={Number of instances (\%)},
          ylabel style={at={(0.03,0.5)}},
          legend style={at={(.97,.9)},anchor=north east},
          legend cell align=left,
          x=0.02604cm, y=0.05cm,
          scale only axis, 
          ymin=0,ymax=104,xmin=0,xmax=250,
          enlargelimits=false,
      ]
      \addplot[mark=none, RedOrange, line width=1pt] table{data/Grid/warm-start/initialSolution.dat};
    \end{axis}
  \end{tikzpicture}
  \qquad
  \begin{tikzpicture}[scale=0.83]
    \begin{axis}[
          xlabel={Optimality gap (\%)},
          ylabel={Number of instances (\%)},
          ylabel style={at={(0.03,0.5)}},
          legend entries={
            \RemoveSpaces{\ref*{formulation:first}},
            \RemoveSpaces{\ref*{formulation:first_splitted}},
            \RemoveSpaces{\ref*{formulation:second}},
            \RemoveSpaces{\ref*{formulation:second_splitted}},
          },
          legend style={at={(.97,.85)},anchor=north east},
          legend cell align=left,
          x=0.0651cm, y=0.08125cm,
          scale only axis, 
          ymin=40,ymax=104,xmin=0,xmax=100,
          enlargelimits=false,
        ]
        \addplot[mark=none, Red, line width=1pt] table{data/Grid/nop/EVSP1.dat};
        \addplot[mark=none, MidnightBlue, line width=1pt] table{data/Grid/nop/EVSP1-S.dat};
        \addplot[mark=none, Green, line width=1pt] table{data/Grid/nop/EVSP2.dat};
        \addplot[mark=none, dashed, Goldenrod, line width=1pt] table{data/Grid/nop/EVSP2-S.dat};
    \end{axis}
  \end{tikzpicture}
  \caption{Percentage of relative optimality gaps between the initial solution value provided by the heuristic and the optimal objective value~(left plot) and percentage of relative optimality gaps between the \mbox{objective} root relaxation and the optimal objective value~(right plot).}
  \label{fig:grid:warm_start}
\end{figure}

\Cref{fig:grid:warm_start} shows performance profiles for the quality of the initial solution obtained by our construction heuristic and the quality of the root relaxation\footnote{The objective value of the root relaxation can be different from that of linear relaxation, since it is improved by tightening the integer linear program given to the solver.}.~This heuristic obtained a solution that has at most~\(2.5\) times the optimal objective value.~On the other plot, we see that the formulations~\RemoveSpaces{\ref{formulation:first_splitted}} and~\RemoveSpaces{\ref{formulation:second_splitted}} stand out as they have better root relaxation.~More precisely,~\(82\%\) of the instances for these formulations achieved values of root relaxation corresponding to the optimal integer solution.~Thus, for these instances, the main difficulty for the solver was to quickly find integer solutions.~It is important to highlight that for~\RemoveSpaces{\ref{formulation:first}} and~\RemoveSpaces{\ref{formulation:second}}, \(33.5\%\) and~\(4\%\) of the instances, respectively, required more than one hour to complete the root relaxation~(the most difficult instance required almost~\(7\) hours).

\subsection{Case study of VAMO Fortaleza} 

For this second performance analysis, we created a benchmark of instances based on real-world data usage of the~VAMO Fortaleza system~\citep{Fortaleza/2018}.~In August~2018, VAMO operated with~\(20\) electric cars,~\mbox{\(5\) BYD e6} and~\mbox{\(15\)~Zhidou EEC}, distributed in~\(12\) stations, but since~December~2019, the Zhidou EEC model has been replaced by the~Renault ZOE.~Since the~\shortname considers a homogeneous fleet, we disregarded the~\mbox{\(5\) BYD e6} and considered~\(15\) Renault ZOE with a~\(52\)~kWh battery, \(3\) hours to full charge the battery, and a discharge rate of~\mbox{\(\lambda = 5.2\) kWh}~(assuming an average speed~\(\vartheta = 30\)~km/h).~Among the~\(12\) stations, two have capacity of~\(3\) and the others have capacity of~\(4\), and all parking spaces are equipped with charging facilities.~The operation period of~VAMO is from~\mbox{6:00} to~\mbox{01:00}~(nineteen hours) per day.

Our second benchmark set is defined considering four different scenarios.~The first scenario is the one closest to the~VAMO Fortaleza setup, where each customer has exactly one demand, and it is divided into groups of~\mbox{\(250, 260, \dots, 340\)} customers.~The second scenario was created assuming that each customer has at most four demands, which are divided into groups of~\mbox{\(30, 40, \dots,\) \(120\)} customers.~By modifying the parameters of the stations used in the second scenario, the third scenario increases the capacity of the stations by varying the number of parking spaces equipped and unequipped with charging facilities, as well as the number of vehicles, as described for the grid instances.~Finally, the last scenario is derived from the third one, in which we changed the capacity of each station to be sufficient for all vehicles in the system.~Each scenario of instances is created as described in the following paragraph (all histograms mentioned are included in the report of city hall).

\noindent \hspace*{-.5mm}\emph{\textbf{VAMO Fortaleza instances}.}~At the first step, all vehicles are randomly distributed among the stations.~Afterward, each demand~\mbox{\((\outgoing{s}, t_{i}, \incoming{s}, t_{j}, \varepsilon) \in \matheuler{D}_{c}\)}, for all~\mbox{\(c \in \matheuler{C}\)}, is defined as follows.~The pick-up and drop-off stations are selected with a certain probability, which is the same as the probability for the origin and destination stations in the city hall report.~The departure time is defined by first randomly selecting one of~\(19\) one-hour intervals of the operating period, where each one-hour interval has a probability proportional to the histogram of the depurate time series, and then it is randomly chosen from the set of time instants~(with a granularity of ﬁve minutes) associated with the selected one-hour interval.~Next, the arrival time is defined as the departure time plus a rental time that is randomly chosen with equal probability from the travel time histogram with a granularity of thirty minutes.~The energy required for the rental period is defined as~\(dis/\vartheta\times\lambda\), where the travel distance~\(dis\) is randomly chosen from the set~\mbox{\(\{5 \text{ km}, 10 \text{ km}, \dots, 70 \text{ km}\}\)} with a certain probability from the travel distance histogram such that~\(dis/\vartheta\) is smaller than the rental time.

\begin{figure}[H]
  \centering
  \resizebox{\textwidth}{!}{%
  \captionsetup[subfigure]{oneside,margin={.5cm,0cm}}
  \begin{subfigure}[t]{.47\textwidth}
    \begin{tikzpicture}[scale=0.83]
      \begin{axis}[
            xlabel={Time (minutes)},
            ylabel={Number of solved instances (\%)},
            ylabel style={at={(0.03,0.5)}},
            legend entries={
              \RemoveSpaces{\ref*{formulation:first}},
              \RemoveSpaces{\ref*{formulation:first_splitted}},
              \RemoveSpaces{\ref*{formulation:second}},
              \RemoveSpaces{\ref*{formulation:second_splitted}},
            },
            legend pos=south east,
            legend cell align=left,
            x=0.651cm, y=0.06cm,
            scale only axis, 
            ymin=0,ymax=104,xmin=0,xmax=10,
            enlargelimits=false,
        ]
        \addplot[mark=none, Red, line width=1pt] table{data/VAMO/Main/performance/EVSP1.dat};
        \addplot[mark=none, MidnightBlue, line width=1pt] table{data/VAMO/Main/performance/EVSP1-S.dat};
        \addplot[mark=none, Green, line width=1pt] table{data/VAMO/Main/performance/EVSP2.dat};
        \addplot[mark=none, Goldenrod, line width=1pt] table{data/VAMO/Main/performance/EVSP2-S.dat};
      \end{axis}
    \end{tikzpicture}
    \caption{Result for Scenario I.}
  \end{subfigure}
  \qquad
  \begin{subfigure}[t]{.47\textwidth}
    \begin{tikzpicture}[scale=0.83]
      \begin{scope}[
          spy using outlines={
              rectangle,
              magnification=3,
              connect spies,
              size=2cm,
              blue,
          },
        ]
        \begin{axis}[
              xlabel={Time (minutes)},
              ylabel={Number of solved instances (\%)},
              ylabel style={at={(0.03,0.5)}},
              legend entries={
                \RemoveSpaces{\ref*{formulation:first}},
                \RemoveSpaces{\ref*{formulation:first_splitted}},
                \RemoveSpaces{\ref*{formulation:second}},
                \RemoveSpaces{\ref*{formulation:second_splitted}},
              },
              legend pos=south east,
              legend cell align=left,
              x=0.1085cm, y=0.06cm,
              scale only axis, 
              ymin=0,ymax=104,xmin=0,xmax=60,
              enlargelimits=false,
          ]
          \addplot[mark=none, Red, line width=1pt] table{data/VAMO/I/performance/EVSP1.dat};
          \addplot[mark=none, MidnightBlue, line width=1pt] table{data/VAMO/I/performance/EVSP1-S.dat};
          \addplot[mark=none, Green, line width=1pt] table{data/VAMO/I/performance/EVSP2.dat};
          \addplot[mark=none, Goldenrod, line width=1pt] table{data/VAMO/I/performance/EVSP2-S.dat};
          \coordinate (spypoint) at (axis cs:44.6,77.6); 
          \coordinate (spyviewer) at (axis cs:10,20); 
          \spy[gray] on (spypoint) in node [fill=white] at (spyviewer);
        \end{axis}
      \end{scope}
    \end{tikzpicture}
    \caption{Result for Scenario II.}
  \end{subfigure}
  }
  \caption{Percentage of instances solved optimally within one hour.~For clarity, the left plot is displayed with an interval of~\(10\) minutes.}
  \label{fig:vamo:part_one}
\end{figure}
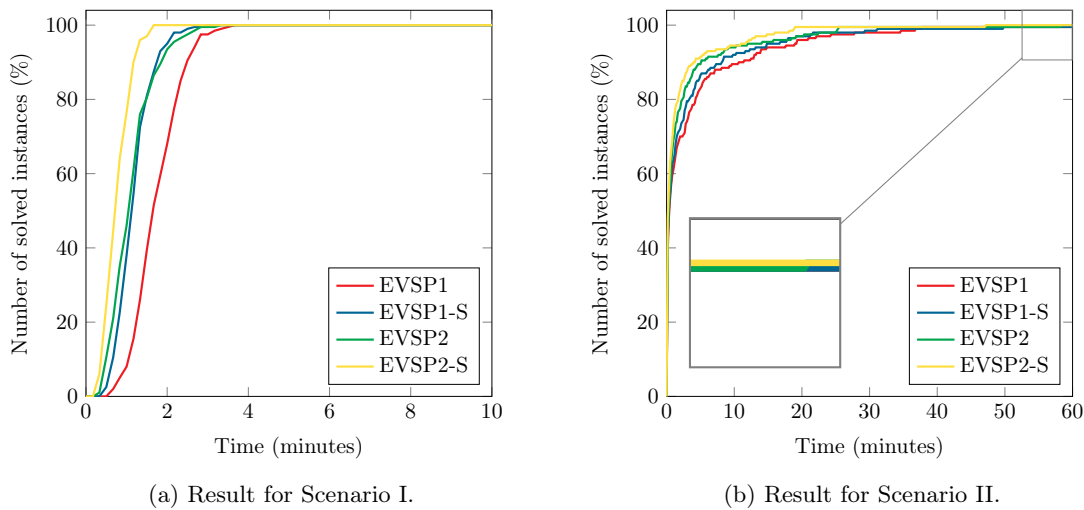

Now we present an overview of the results.~\Cref{fig:vamo:part_one} shows the percentage of instances solved for the first two scenarios.~For the scenario where each customer has only one demand~(left plot), all instances were solved very quickly, in less than~\(4\) minutes.~For the scenario where customers can have one or more demands~(right plot), only~\RemoveSpaces{\ref*{formulation:first}} and~\RemoveSpaces{\ref*{formulation:first_splitted}} did not achieve an optimal solution for all instances.~More precisely, just for one instance they did not achieve an optimal solution, delivering the lower and upper bounds with a gap of almost~\(4\%\).~We observed that for these two scenarios, as all parking spaces have charging facilities, all formulations have an advantage since the variables~\(y^{v}_{e}\) are removed from their integer linear programs.

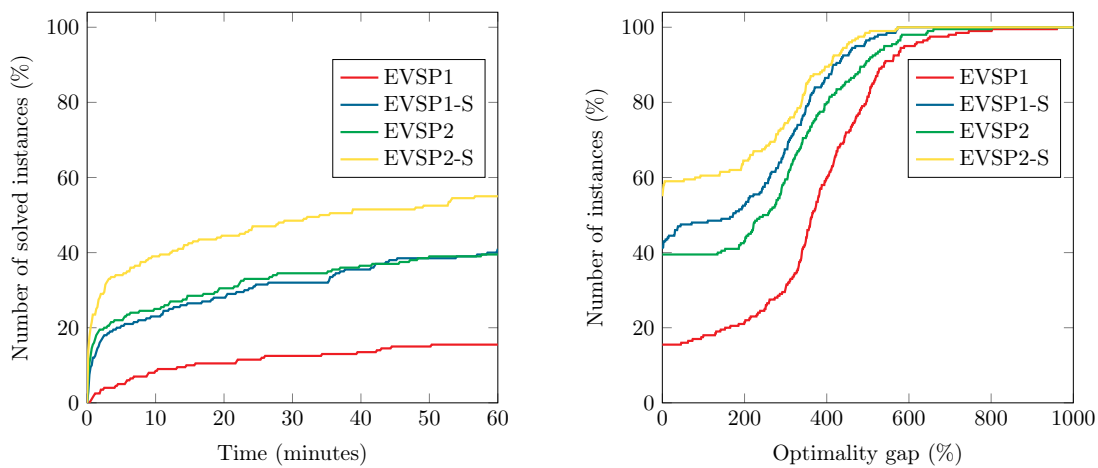
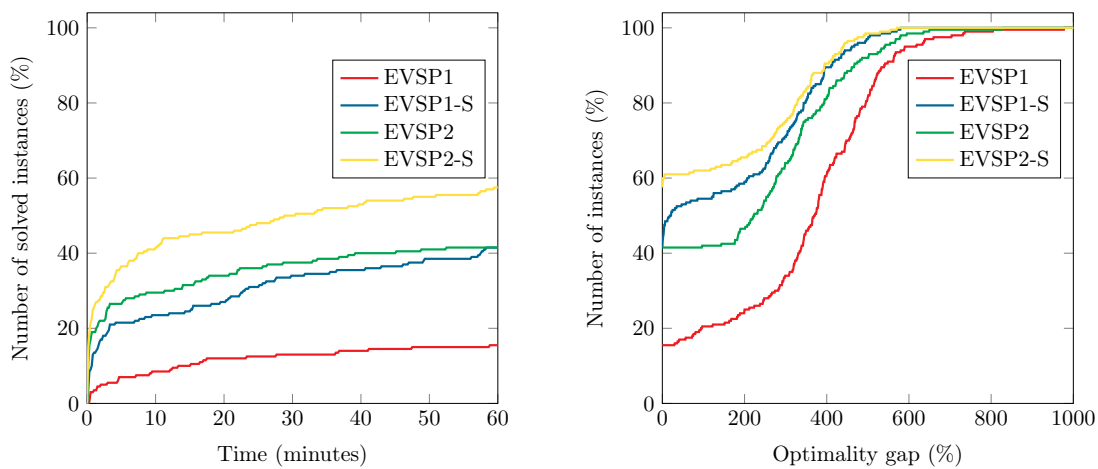
\begin{figure}[H]
  \centering
  \begin{subfigure}[t]{\textwidth}
    \begin{tikzpicture}[scale=0.83]
      \begin{axis}[
            xlabel={Time (minutes)},
            ylabel={Number of solved instances (\%)},
            ylabel style={at={(0.03,0.5)}},
            legend entries={
              \RemoveSpaces{\ref*{formulation:first}},
              \RemoveSpaces{\ref*{formulation:first_splitted}},
              \RemoveSpaces{\ref*{formulation:second}},
              \RemoveSpaces{\ref*{formulation:second_splitted}},
            },
            legend style={at={(.97,.88)},anchor=north east},
            legend cell align=left,
            x=0.1085cm, y=0.06cm,
            scale only axis, 
            ymin=0,ymax=104,xmin=0,xmax=60,
            enlargelimits=false,
        ]
        \addplot[mark=none, Red, line width=1pt] table{data/VAMO/II/performance/EVSP1.dat};
        \addplot[mark=none, MidnightBlue, line width=1pt] table{data/VAMO/II/performance/EVSP1-S.dat};
        \addplot[mark=none, Green, line width=1pt] table{data/VAMO/II/performance/EVSP2.dat};
        \addplot[mark=none, Goldenrod, line width=1pt] table{data/VAMO/II/performance/EVSP2-S.dat};
      \end{axis}
    \end{tikzpicture}
    \qquad
    \begin{tikzpicture}[scale=0.83]
      \begin{axis}[
            xlabel={Optimality gap (\%)},
            ylabel={Number of instances (\%)},
            ylabel style={at={(0.03,0.5)}},
            legend entries={
              \RemoveSpaces{\ref*{formulation:first}},
              \RemoveSpaces{\ref*{formulation:first_splitted}},
              \RemoveSpaces{\ref*{formulation:second}},
              \RemoveSpaces{\ref*{formulation:second_splitted}},
            },
            legend style={at={(.97,.88)},anchor=north east},
            legend cell align=left,
            x=0.00651cm, y=0.06cm,
            scale only axis, 
            ymin=0,ymax=104,xmin=0,xmax=1000,
            enlargelimits=false,
          ]
          \addplot[mark=none, Red, line width=1pt] table{data/VAMO/II/gap/EVSP1.dat};
          \addplot[mark=none, MidnightBlue, line width=1pt] table{data/VAMO/II/gap/EVSP1-S.dat};
          \addplot[mark=none, Green, line width=1pt] table{data/VAMO/II/gap/EVSP2.dat};
          \addplot[mark=none, Goldenrod, line width=1pt] table{data/VAMO/II/gap/EVSP2-S.dat};
      \end{axis}
    \end{tikzpicture}
    \caption{Results for Scenario III --- Variable number of parking spaces and vehicles.}
  \end{subfigure}

  \medskip\medskip\medskip
  \begin{subfigure}[t]{\textwidth}
    \begin{tikzpicture}[scale=0.83]
      \begin{axis}[
            xlabel={Time (minutes)},
            ylabel={Number of solved instances (\%)},
            ylabel style={at={(0.03,0.5)}},
            legend entries={
              \RemoveSpaces{\ref*{formulation:first}},
              \RemoveSpaces{\ref*{formulation:first_splitted}},
              \RemoveSpaces{\ref*{formulation:second}},
              \RemoveSpaces{\ref*{formulation:second_splitted}},
            },
            legend style={at={(.97,.88)},anchor=north east},
            legend cell align=left,
            x=0.1085cm, y=0.06cm,
            scale only axis, 
            ymin=0,ymax=104,xmin=0,xmax=60,
            enlargelimits=false,
        ]
        \addplot[mark=none, Red, line width=1pt] table{data/VAMO/III/performance/EVSP1.dat};
        \addplot[mark=none, MidnightBlue, line width=1pt] table{data/VAMO/III/performance/EVSP1-S.dat};
        \addplot[mark=none, Green, line width=1pt] table{data/VAMO/III/performance/EVSP2.dat};
        \addplot[mark=none, Goldenrod, line width=1pt] table{data/VAMO/III/performance/EVSP2-S.dat};
      \end{axis}
    \end{tikzpicture}
    \qquad
    \begin{tikzpicture}[scale=0.83]
      \begin{axis}[
            xlabel={Optimality gap (\%)},
            ylabel={Number of instances (\%)},
            ylabel style={at={(0.03,0.5)}},
            legend entries={
              \RemoveSpaces{\ref*{formulation:first}},
              \RemoveSpaces{\ref*{formulation:first_splitted}},
              \RemoveSpaces{\ref*{formulation:second}},
              \RemoveSpaces{\ref*{formulation:second_splitted}},
            },
            legend style={at={(.97,.88)},anchor=north east},
            legend cell align=left,
            x=0.00651cm, y=0.06cm,
            scale only axis, 
            ymin=0,ymax=104,xmin=0,xmax=1000,
            enlargelimits=false,
          ]
          \addplot[mark=none, Red, line width=1pt] table{data/VAMO/III/gap/EVSP1.dat};
          \addplot[mark=none, MidnightBlue, line width=1pt] table{data/VAMO/III/gap/EVSP1-S.dat};
          \addplot[mark=none, Green, line width=1pt] table{data/VAMO/III/gap/EVSP2.dat};
          \addplot[mark=none, Goldenrod, line width=1pt] table{data/VAMO/III/gap/EVSP2-S.dat};
      \end{axis}
    \end{tikzpicture}
    \caption{Results for Scenario IV --- Each station has capacity to cover all vehicles in the system.}
  \end{subfigure}
  \caption{Percentage of optimally solved instances and relative optimality gaps within the one-hour time limit.~Gap values were calculated as~\mbox{(\hspace*{-.1mm}(\hspace*{-.1mm}\textsf{UB} - \textsf{LB})/\textsf{LB})\(\times 100\))}, where \textsf{UB} and~\textsf{LB} are the upper and lower bound values achieved by the solver, respectively.}
  \label{fig:vamo:part_two}
\end{figure}

In the following,~\Cref{fig:vamo:part_two} shows the percentage of solved instances and the remaining optimality gaps for the last two scenarios.~If we compare these scenarios, we can see that the plots are similar to each other.~In other words, changing the capacity of each station so that it is sufficient to cover all vehicles in the system has relatively little impact on the optimization process.

\section{Conclusion}
\label{sec:conclusions}

In this paper, we introduced the electric vehicle sharing problem~(\shortname) which is motivated by planning and operating one-way electric car-sharing systems.~We proposed four mixed-integer linear programming formulations for this problem, based on space-time networks.~The first two formulations are based on homogeneous networks, which means that the same network is considered for all vehicles in the system.~These formulations are referred to as~\RemoveSpaces{\ref{formulation:first}} and~\RemoveSpaces{\ref{formulation:first_splitted}}, respectively.~The last two formulations are based on heterogeneous networks, which means that a specific network is considered for each vehicle in the system.~Roughly speaking, it involves several steps of preprocessing to reduce the size of the networks.~These formulations are referred to as~\RemoveSpaces{\ref{formulation:second}} and~\RemoveSpaces{\ref{formulation:second_splitted}}, respectively.~We showed theoretical results in terms of the strength of their linear programming relaxation.

We also presented a complexity analysis showing that the~\shortname is an~\mbox{NP-hard} problem even if there is one station and each customer has only one driving demand.

Last, we conducted extensive computational experiments to investigate the performance of our formulations.~For this purpose, we created two benchmark sets.~The first set was created based on grid instances.~On this set, a computational study was performed to analyze the characteristics of the integer linear programs.~In particular, for the formulation~\RemoveSpaces{\ref{formulation:first_splitted}}, the presolve phase of the commercial solver reduced the number of binary variables and constraints by~\(61.69\%\) and~\(69.83\%\) on average, respectively.~However, for these instances, the size of the integer linear program for the formulation~\RemoveSpaces{\ref{formulation:second_splitted}} is smaller than for~\RemoveSpaces{\ref{formulation:first_splitted}} after the presolve phase.~Since the formulations~\RemoveSpaces{\ref{formulation:first_splitted}} and~\RemoveSpaces{\ref{formulation:second_splitted}} have basically the same constraint sets, our ideas are promising as they can reduce the integer linear program better than the commercial solver used.

Another set of instances is based on data from the~VAMO Fortaleza system, the first public one-way electric car-sharing system in Brazil, located in the city of~Fortaleza.~On this set, a computational study was carried out to analyze the efficiency in a real~(small) system and the impact of different scenarios on the overall performance.~For the simplest scenario, where all customers have only one demand, the formulations~\RemoveSpaces{\ref{formulation:first_splitted}}, \RemoveSpaces{\ref{formulation:second}}, \RemoveSpaces{\ref{formulation:second_splitted}}, and even our worst formulation in general~\RemoveSpaces{\ref{formulation:first}} proved to be very effective and solved to optimality all the instances.~On the other hand, for the two most difficult scenarios, our best formulation, that is, the~\RemoveSpaces{\ref{formulation:second_splitted}}, solved to optimality~\(55\%\) and~\(57.5\%\) of the instances, respectively, within the one-hour time limit.

Interesting directions for future research can be derived from our study.~From a computational point of view, we could consider developing a more compact formulation~(in terms of the number of variables) and decomposition methods, such as the Benders decomposition and the Dantzig-Wolfe decomposition, which are likely to produce an exact algorithm that performs significantly better than our current approaches.~Metaheuristic methods could also be used to solve large instances in reasonable computation time.

Apart from these, relevant research directions could consider the relocation of vehicles by operators and the incorporation of demands prioritization, which could involve relaxing the assumption that either all customer demands are completely fulfilled or the customer does not use the system at all.

\section*{Acknowledgments}


This study was financed in part by the~Coordination for the Improvement of Higher Education Personnel --~Brasil (CAPES) --~Finance Code 001.~Supported by grant number \mbox{2015/11937-9}, São Paulo Research Foundation~(FAPESP); and grant numbers \mbox{425340/2016-3}, \mbox{314384/2018-9}, \mbox{425806/2018-9}, \mbox{435520/2018-0} and~\mbox{311039/2020-0}, Brazilian National Council for Scientific and Technological Development~(CNPq).

\bibliographystyle{abbrvnat}
\bibliography{references}

\appendix

\section*{Appendix.\quad Proofs of propositions}
\label{appendix}

In this appendix we give the proofs of the remaining propositions in this paper.~For the sake of clarity, we repeat their statements.

\noindent \textbf{\Cref*{proposition:EVSP1-SL_is_stronger}.}~\textit{%
  There are instances for which~\textnormal{opt(\RemoveSpaces{\ref*{formulation:first_splitted}$_{\textnormal{L}}$})} is smaller than~\textnormal{opt(\RemoveSpaces{\ref*{formulation:first}$_{\textnormal{L}}$})}.%
}

{\noindent {\bf Proof.}~It is sufficient to give an example here.~Consider the following small-size instance with two stations, denoted by~\(0\) and~\(1\), where each station having one parking space equipped with a charging facility and another unequipped.~At each station, there is initially one vehicle with a battery capacity of~\(600\) \mbox{watt-minutes~(Wm), with} a maximum charging time of one hour.~In this instance, there are three customers, each with one of the following sets of demands:~\mbox{\(\{(1, 745, 0, 861, 465),\, (0, 900, 1, 1011, 530)\}\)}, \mbox{\(\{(0, 495, 1, 591, 475),\, (1, 760, 0, 871, 455)\}\),} \mbox{and \(\{(0, 895, 1, 991, 455),\, (1, 1120, 0, 1241, 485)\}\)}. For this instance,~\Cref{fig:EVSP1,fig:EVSP1-S} illustrate an optimal solution for the formulations.~In these illustrations, for clarity, some nodes are removed and the connecting arcs between the removed nodes are combined into one.

\begin{figure}[H]
  \centering
  \resizebox{\textwidth}{!}{%
    \input{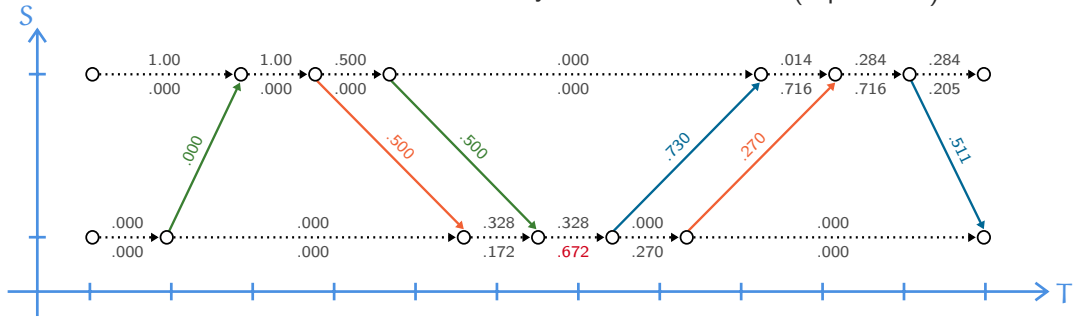}
  }
  \caption{Illustration of an optimal fractional solution for~\RemoveSpaces{\ref*{formulation:first}$_{\textnormal{L}}$}.~Here, there is one simplified network for each vehicle.~Each demand arc is labeled with the value assigned to the variable~\(x^{v}_{e}\), whereas each connecting arc is labeled above and below with the value assigned to the variables~\(y^{v}_{e}\) and~\(z^{v}_{e}\), respectively.~In addition, note that the capacity of the charging facilities is exceeded when the red values on connecting arcs are added.}
  \label{fig:EVSP1}
\end{figure}

As stated, an important drawback of~\RemoveSpaces{\ref*{formulation:first}$_{\textnormal{L}}$} is that it can have connecting arcs whose corresponding flows exceed the number of parking spaces equipped or unequipped with charging facilities.

\begin{figure}[H]
  \centering
  \resizebox{\textwidth}{!}{%
    \input{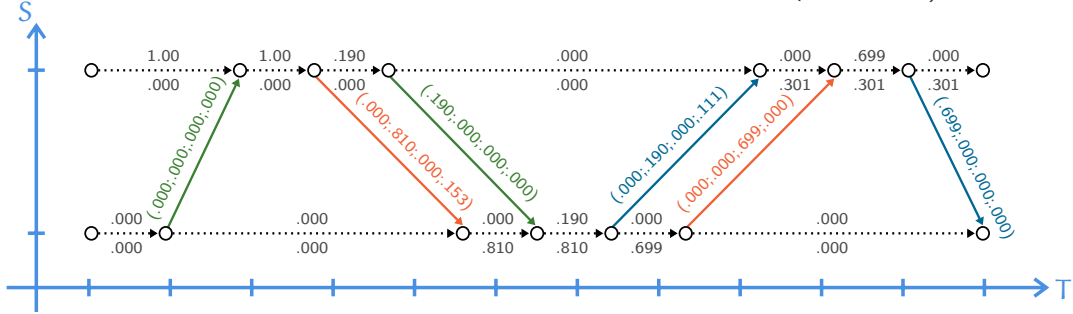}
  }
  \caption{Illustration of an optimal fractional solution for~\RemoveSpaces{\ref*{formulation:first_splitted}$_{\textnormal{L}}$}.~Here, there is one simplified network for each vehicle.~Each demand arc is labeled with the value assigned to the variables~\(x^{v}_{e\hspace*{-.2mm},\hspace*{-.2mm}k}\), for all~\mbox{\(k \in \matheuler{K}\)}, whereas each  connecting arc is labeled above and below with the value assigned to the variables~\(y^{v}_{e}\) and~\(z^{v}_{e}\), respectively.}
  \label{fig:EVSP1-S}
\end{figure}

\vspace*{-.25cm}
Since the corresponding optimal objective value for~\RemoveSpaces{\ref*{formulation:first}$_{\textnormal{L}}$} and~\RemoveSpaces{\ref*{formulation:first_splitted}$_{\textnormal{L}}$} are not the same for this example, we obtain the result.}~\hfill \(\square\)

\medskip\medskip

\noindent \textbf{\Cref*{proposition:EVSP2-SL_is_stronger}.}~\textit{%
  There are instances for which~{\normalfont opt(\RemoveSpaces{\ref*{formulation:second_splitted}$_{\textnormal{L}}$})} is smaller than~{\normalfont opt(\RemoveSpaces{\ref*{formulation:second}$_{\textnormal{L}}$})}.%
}

{\noindent {\bf Proof.}~Consider the same example used in~\Cref{proposition:EVSP1-SL_is_stronger}.~Since the corresponding optimal \mbox{objective} value for~\RemoveSpaces{\ref*{formulation:second_splitted}}$_{\textnormal{L}}$ is~\(625.99\) and the optimal objective value for~\RemoveSpaces{\ref*{formulation:second}}$_{\textnormal{L}}$ is~\(630.10\), the result follows.}~\hfill \(\square\)

\end{document}